\documentclass{amsart}
\usepackage{amsmath,amssymb,MnSymbol}
\usepackage{graphicx,color}
\usepackage[mathscr]{euscript}
\usepackage{paralist}
\usepackage{subcaption,caption}

\newcommand{\scp}[2]{{\left\langle {#1}\, , \, {#2}\right\rangle}}
\newcommand{\lform}[2]{{\left( {#1}\, \left\vert \vphantom{#1}\,  {#2}\right.\right)}}
\newcommand{\Diff}{{\mathit{Diff}\!}}
\newcommand{\Q}{{\mathscr{Q}}}
\newcommand{\G}{{\mathscr{G}}}
\newcommand{\K}{{\mathscr{K}}}

\renewcommand{\H}{{\mathscr{H}}}
\newcommand{\M}{{\mathscr{M}}}
\newcommand{\mR}{{\mathbb R}}
\newcommand{\id}{\mathrm{id}}
\renewcommand{\ker}{\mathit{Null}}
\newcommand{\normi}[1]{{\left\lWavy \vphantom{|} #1\right\rWavy}}
\graphicspath{{./}}
\usepackage[numbers]{natbib}

\title{Hybrid Riemannian Metrics for Diffeomorphic Shape Registration}

\author{Laurent Younes}
\address{Department of Applied Mathematics and Statistics\\
 and Center for Imaging Science,\\
  Johns Hopkins University\\}

\begin{document}
\maketitle

\begin{abstract}
We consider the results of combining two approaches developed for the design of Riemannian metrics on curves and surfaces, namely parametriza\-tion-invariant metrics of the Sobolev type on spaces of immersions, and metrics derived through Riemannian submersions from right-invariant Sobolev metrics on groups of diffeomorphisms (the latter leading to the ``large deformation diffeomorphic metric mapping'' framework). We show that this quite simple approach inherits the advantages of both methods, both on the theoretical and experimental levels, and provide additional flexibility and modeling power, especially when dealing with complex configurations of shapes. Experimental results illustrating the method are provided for curve and surface registration.
\end{abstract}

\section{Introduction}
\subsection{Shape Registration: Basic Principles}

We consider  a ``shape space'' (denoted $\M$) in which objects are subject to free-form deformations, so that  a group action 
$
(\varphi, q)\in \Diff \times \M \mapsto \varphi\cdot q \in \M
$
is defined on $\M$(where $\Diff$ refers to the diffeomorphism group). This concept usually comes with additional conditions on the structure of the space. Here, following  \cite{arguillere2015abstract}, we will assume that $\M$ is an open subset of a Banach space $\Q$. One also often considers $\M$ quotiented by other group actions (such as Euclidean transformations, or reparametrization). Even though we will not formally consider such quotient spaces,  such invariance will often be a direct consequence of the models we will discuss.

%The term LDDMM represents a collection of shape analysis algorithms,  in which shapes  are registered to a template in order to provide local coordinates. These coordinates can then be used as a multivariate representation, allowing, for example, for statistical  analyses.

One can interpret registration methods within the following  framework. 
%Consider a ``shape space'', $\mathscr Q$, whose elements of $q$ are typically point sets, curves, surfaces, or images. 
%Recall that a diffeomorphism of an open subset $\Omega \subset \mathbb R^n$ is a $C^1$ transformation $\varphi: \Omega\to\Omega$,  which is invertible with a $C^1$ inverse. We will denote  $\Diff$ as the set of diffeomorphisms of $\mathbb R^n$, which forms a group, and assume that this group acts on $\Q$ with an action denoted $(\phi, q) \mapsto \phi\cdot q$.
Given a ``template'' $q_0\in \M$, a registration method can be seen as a transformation that takes a shape $q$ as input and returns as output a diffeomorphism $\varphi$ such that $\varphi\cdot q_0 = \text{(or }\simeq \text{) } q$, therefore providing a mapping, $\Psi_{q_0}: \M\to \Diff$, representing each shape by a diffeomorphism. 
%Therefore, it provides a representation of the shape space within the diffeomorphism group.
%, this representation being relative to a fixed template, who is itself mapped to the identity. 
One of the main advantages of such constructions is that it is much easier to define characteristic features of diffeomorphisms than of shapes,  considered, for example, as subsets of $\mR^d$. Voxel-based, or surface-based morphometric methods have exploited this by introducing deformation markers,  often deduced from the Jacobian of the estimated diffeomorphism \cite{af00,pantazis2004statistical}.

%It is important to realize that the set of diffeomorphisms that may represent a shape does not cover in general the full diffeomorphism group.
 When working with shape spaces of landmarks, images, curves or surfaces, there exist, for each given shape $q$,  either zero or an infinity of transformations such that $\varphi\cdot q_0 = q$. They can all be deduced from each other via composition on the right by diffeomorphisms that leave $q_0$ invariant, i.e., elements of the stabilizer
 $
 \mathcal S_{q_0} = \{\psi\in\Diff: \psi\cdot q_0 = q_0\}.
 $ 
 ``Good'' registration algorithms generally pick one such transformation that maximizes a regularity or optimality criterion that the algorithm implements. Understanding the structure of the space $\Psi_{q_0}(\Q)$ has many advantages, because it may lead to (locally) one-to-one shape representations. Among others methods \cite{christensen1996deformable,miller2001group,droske2004variational,trouve2005metamorphoses,styner2006framework,ashburner2007fast,hernandez2008comparing,vercauteren2009diffeomorphic,klein2009evaluation,holm2009euler,ashburner2011diffeomorphic}, the large deformation diffeomorphic metric mapping framework (LDDMM) includes a family of registration algorithms \cite{jos97,joshi2000Landmark,avants2004geodesic,beg2005computing,cao2005large,glaunes2008large,ceritoglu2009multi,risser2011simultaneous}, adapted to various shape modalities,  that rely on a sub-Riemannian setup of the diffeomorphism group and of the shape space. This setup is a special case of the model used in this paper, that we now describe.

 We will assume a sub-Riemannian structure on $\Diff$ and consider conditions under which this action can be transferred into a sub-Riemannian structure on $\M$. This framework will include the LDDMM construction as a special case, and the other metrics that will be used in this paper.  The following notation and assumptions will be used throughout this paper. 
 
 We will only consider diffeomorphisms that tend to the identity at infinity, i.e., $\varphi = \id + u$ such that $u$ and $du$ tend to 0 at infinity. $\Diff$ will denote the space of such diffeomorphisms, and $u\mapsto \id+u$ provides a trivial chart of $\Diff$ as a Banach manifold, when defined over the space $C_0^1(\mR^d, \mR^d)$ of $u$'s that tend to 0 at infinity (with the supremum norm: $\|u\|_{1, \infty} = \|u\|_\infty + \|du\|_\infty$). We will assume that the action $(\varphi, q) \mapsto \varphi\cdot q$ is $C^1$ from $\Diff\times \M \to \M$. We let $\pi_q(\varphi) = \varphi\cdot q$ and $\xi_q = d\pi_q(\id)$, so that the infinitesimal action is given by $v\cdot q =  \xi_qv$.

 Let $V$ be a Hilbert space continuously embedded in $C_0^1(\mR^d, \mR^d)$. Consider the sub-Riemannian structure on $\Diff$ associated with the distribution $V_\varphi = \{v\circ \varphi: v\in V\}$. We assume that $V_\varphi$ is equipped with a Hilbert structure, with norm $\|\cdot \|_\varphi$ such that, for all $v\in V$,  $\|v\circ \varphi\|_\varphi \geq c \|v\|_V$ for some $c>0$. We will also denote $\|v\|_{V, \varphi} = \|v\circ \varphi\|_\varphi$ for $v\in V$. A path $(\varphi(t), t\in [0,1])$ in $\Diff$ is admissible if $\partial_t\varphi \in V_\varphi$ for (almost all) $t$ and 
 \[
 \int_0^1 \|\partial_t\varphi \|^2_{\varphi(t)} \,dt <\infty.
 \]
 (Paths tangent to $V_{\varphi(t)}$  at all times are often called horizontal in the sub-Riemannian literature. We will not use this term here because of the horizontal spaces that we define just below.)
 
 Let $\Diff_0$ be the subgroup of $\Diff$ containing all the elements that are reachable from the identity with an admissible path. Fix $q_0\in \M$ and let $\M_0= \{\pi_{q_0}(\varphi):\varphi\in\Diff_0\}$. For $\varphi\in \Diff_0$ and $q = \pi_{q_0}(\varphi)$, consider the space orthogonal to $\ker(\xi_q)$ for $\scp{\cdot}{\cdot}_{V, \varphi}$, denoted
$H_\varphi$.
Assume  that the Hilbert space isometry
\begin{equation}
\label{eq:sub.sub}
(H_\varphi, \|\cdot\|_{V, \varphi}) \sim (H_\psi, \|\cdot\|_{V, \psi})
\end{equation}
holds whenever  $\pi_{q_0}(\varphi)=\pi_{q_0}(\psi)$.
Then, for $\pi_{q_0}(\varphi) = q$, one can define the space $\H_q = \xi_qH_\varphi = \{\xi_qv: v\in H_\varphi\}$ with 
\[
\|\xi_q v\|_q =  \|v\|_{V,\varphi}
\]
and this definition does not depend on which $\varphi\in \pi_{q_0}^{-1}(q)$ is chosen.
The distribution $\H_q$ then provides  a sub-Riemannian structure on $\M_0$. 

The space $\H_\varphi = \{v\circ \varphi: v\in H_\varphi\}$ is the horizontal space at $\varphi$ for the mapping $\pi_{q_0}$ and the statement that these spaces are isometric adapts the conditions for $\pi_{q_0}$ to be a Riemannian submersion to this sub-Riemannian setting. In the LDDMM framework, \eqref{eq:sub.sub} is ensured by defining $\|v\circ \varphi\|_\varphi = \|v\|_V$ for all $\varphi$ and $v\in V$, so that the original metric is right-invariant. We will below consider computationally feasible situations in which the latter condition is relaxed with \eqref{eq:sub.sub} still holding.\\

From a practical viewpoint,  LDDMM  can be expressed as an optimal control problem. One can indeed describe the search for a geodesic between the template $q_0$ and a shape $q_1\in \M$ as the minimization, over all time-dependent vector fields $v:[0,1]\to V$, of
\begin{equation}
\label{eq:exact.1}
\frac12 \int_0^1 \|v(t)\|^2_V\,dt
\end{equation}
subject to $\partial_t q(t) = v(t)\cdot q(t)$,  $q(0) = q_0$, and $q(1) = q_1$.
LDDMM uses a relaxation of the last constraint, minimizing
\begin{equation}
\label{eq:lddmm}
\frac12 \int_0^1 \|v(t)\|^2_V\,dt + D(q(1), q_1)  
\end{equation}
subject to $\partial_t q(t) = v(t)\cdot q(t)$ and $q(0) = q_0$,
 where $D(q,\tilde q)$ is some properly defined discrepancy measure between $q$ and $\tilde q$. Invariance is often ensured by considering functions such that $D(q,\tilde q)=0$ if $q$ and $\tilde q$ are related by a transformation for which the invariance is searched.
 
\begin{figure}
\centering
\includegraphics[width=0.85\textwidth]{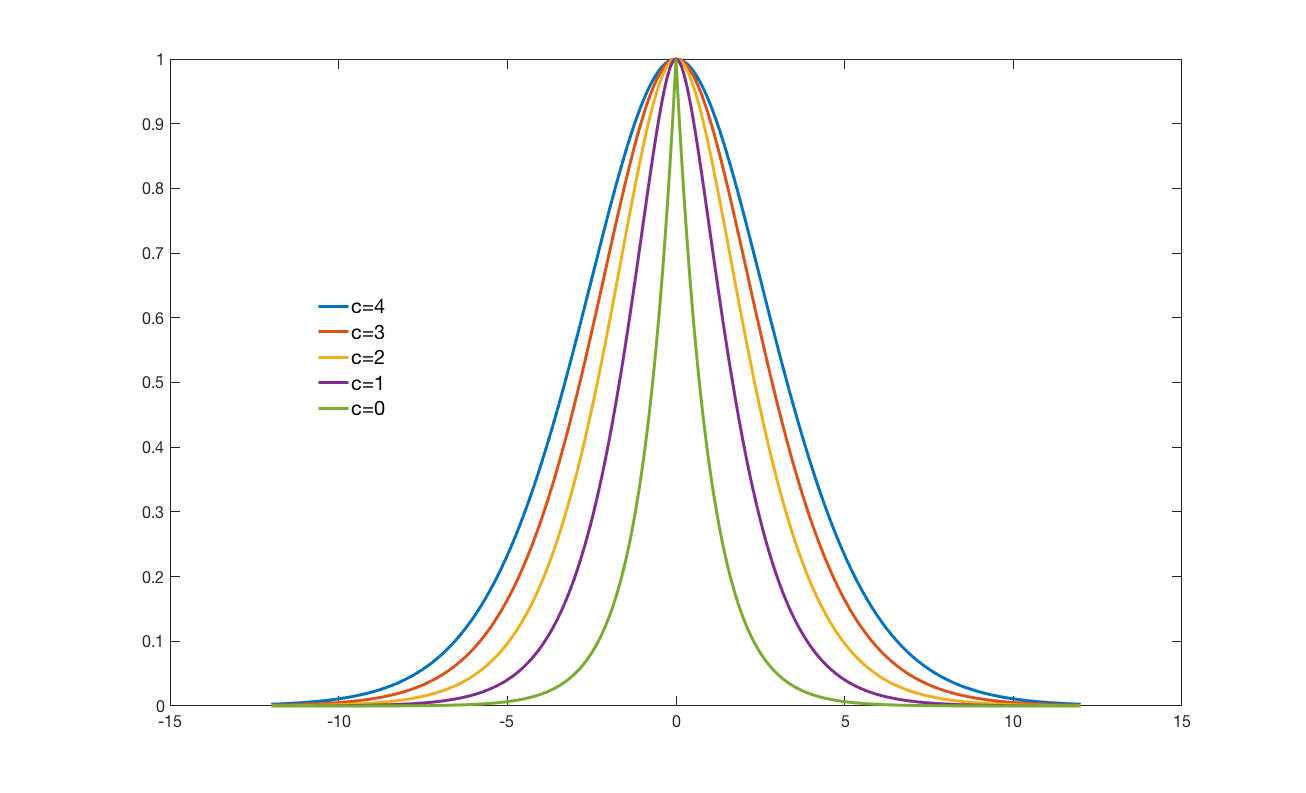}
\caption{\label{fig:kernels} Shapes of the kernels given by \eqref{eq:abel.kernel} for various values of $c$.}
\end{figure}

Because of the embedding assumption, the norm on $V$ is associated to a reproducing kernel, which is a matrix-valued function:
\[
\begin{aligned}
K:&& \mR^d\times\mR^d & \to && \mathcal M_d(\mR^d)\\
&& (x,y) & \mapsto && K(x,y)
\end{aligned}
\]
such that, for all $x,a\in \mR^d$,  $K(\cdot,x)a: y \mapsto K(y,x)a$ belongs to $V$, and for all $v\in V$,
\[
a^T v(x) = \scp{v}{K(\cdot, x)a}_V.
\]
This kernel and the norm on $V$ are generally chosen to be translation invariant, taking the form
\begin{equation}
\label{eq:gamma}
K(x,y) = \Gamma\left(\frac{x-y}{a}\right)
\end{equation}
where $\Gamma$ is a positive definite (matrix-valued) function. The extra parameter, $a$, can be interpreted as a scale parameter, that can be tuned to modulate the locality of the deformations. It essentially modulates the long range effect of the motion of  a single particle in space. For example, the kernel associated with
\[
\|v\|_V^2 = \int_{\mR^d}\left|(\mathrm{Id} - a^2\Delta)^{m/2} v\right|^2\, dx,
\]
where $c=m-(d+1)/2$ is a positive integer, is given by \eqref{eq:gamma} with
\begin{equation}
\label{eq:abel.kernel}
\Gamma(x) = P_c(|x|)e^{-|x|},
\end{equation}
where $P_c$ is a  {\em reverse Bessel polynomial} of degree $c$ (see \cite{krall1949new}), normalized so that $P(0)=1$. The associated kernel $K$ decays to 0 at infinity, at a speed which is modulated by the scale constant $a$. The shape of the function $\Gamma$ for $a=1$ and $c=0,1, \ldots,4$ is provided in Figure \ref{fig:kernels} ($\Gamma$ has $c$ continuous derivatives at $x=0$ and is $C^\infty$ everywhere else). We used $c=3$ in our experiments. For this kernel, the half-range (value of $|x|$ for which $\Gamma(x) = 1/2$) is given by $2.85a$. 

The rather simple formulation (ignoring  numerical issues) provides a horizontal geodesic in $\Diff_0$, given by the flow associated to an optimal $v$, i.e., the solution of $\partial_t \varphi(t) = v(t)\circ \varphi(t)$ with initial condition $\varphi(0) = \mathit{id}$. 
%Because we have a Riemannian submersion, this geodesic is, in addition, horizontal, meaning that $v(t)\circ \varphi(t) \in H_{\varphi(t)}$ at all times. This implies that $\Psi_{q_0}(\Q)$ can be described as the set of diffeomorphisms that can be reached by an horizontal geodesic starting at $q_0$, and since geodesics are uniquely specified by their first-order initial conditions, an element of $\Psi_{q_0}(\Q)$ can be described by the initial velocity $v(0)\in H_{q_0}$ of a geodesic that reaches it from $q_0$. One can then use element of $H_{q_0}$ as coordinates parametrizing the shape space, obtaining in that way a representation isometric to the exponential chart at $q_0$ in the Riemannian manifold $\Q$. For $v\in H_{q_0}$, the covariant form $v^\flat: w \mapsto \scp{v}{w}_V$, provides an equivalent formulation which turns out to be much more parsimonious than $v$. The so-called {\em initial momentum}, $v^\flat(0)$,  has become a preferred way of representing shapes though LDDMM. For landmark shapes spaces, it can be parametrized with a finite number of $d$-dimensional vectors associated with each landmark (while $v(0)$ is a vector field defined on the whole $\mR^d$). For images, $v^\flat(0)$ can be associated to a vector field that is normal to level sets of the image, etc.

%\subsection{Our methods and innovations} \label{Sec2.3}
%\subsection*{Relaxing the Right Invariance Condition}
\section{Hybrid LDDMM}
Starting with a right-invariant metric on $\Diff_0$ ensures that all tangent spaces are isometric to the tangent space at the identity, through the right-translation map $v \mapsto v \circ \varphi$. This property is much stronger than what is needed to ensure \eqref{eq:sub.sub}, which only requires that horizontal spaces within the same fibers to be isometric. Obviously, right-invariance brings additional properties to the Riemannian structure on $\Diff_0$, making it, in particular, independent of the choice of the template, $q_0$, and ensuring that its geodesics satisfy strong conservation laws \cite{holm1998euler,younes2009evolutions,miller2006geodesic,younes2007jacobi}. 
On the other hand, it prevents the metric from taking into account  shape-dependent properties, related, for example to the geometry of the considered curves or surfaces. 
%The way LDDMM defines optimal transformations of, say, curves in 2D is by computing a global deformation of the whole space, typically considered as a homogeneous medium with minimal energy subject to 

%To ensure that the submersion $\pi$ is Riemannian, it suffices that the metric at $\psi\in \Diff_0$ takes the form
%$
%\|h\|_\psi = \| h\circ \psi^{-1}\|_{\psi\cdot q_0}, 
%$ 
%where $\|\cdot\|_q$ is a family of norms on $V$ indexed by $q\in \Q$. Geodesics in $\Q$ can now be associated with the minimization of 
%\begin{equation}
%\label{eq:exact.2}
%\frac12 \int_0^1 \|v(t)\|^2_{q(t)}\,dt 
%\end{equation}
%subject to $\partial_t q(t) = v(t)\circ q(t)$,  $q(0) = q_0$, and $q(1) = q_1$.
%
% 
%
Allowing for less restrictive invariance will allow us to characterize a much larger variety of features compared to those associated with plain LDDMM. This is done in the following examples over spaces of curves and surfaces.

%In the following examples, we assume that a fixed Hilbert norm $\|\cdot\|_V$ is chosen on $V$, and that this norm ensures sufficient smoothness properties for the flow defined by $\partial_t \varphi(t) = v(t)\circ \varphi(t)$ to be well defined as soon as $\int_0^1 \|v(t)\|^2_V\,dt<\infty$. All shape-dependent norms $\|\cdot\|_q$ will be assumed to be continuous, ensuring that the condition $\int_0^1 \|v(t)\|^2_{q(t)}\,dt<\infty$ also implies the well-posedness of the flow. For $q\in\M$, we will denote by $\xi_q: V \to \Q$ the operator $v\mapsto v\cdot q$.  

\subsection{Curves}

\subsubsection{Hybrid Norms}
Consider  the situation in which the objects of interest are curves, in two or three dimensions. In this setting, we can take advantage of the collection of metrics that have been introduced for shape spaces of embedded curves $q: M \to \mR^2$, where $M$ is either the unit interval or the unit disk. We let $\Q = C^r(M, \mR^d)$ (normed by the supremum norm over all derivatives of order $r$ or less), for some $r>1$, and $\M$ be the set of $C^r$ embeddings. The action being $\varphi\cdot q = \varphi\circ q$, we have $\xi_q v = v\circ q$, and we will assume that $V$ is continuously embedded in $C^p_0(\mR^d, \mR^d)$ with $p>r$, which ensures that $(v,q)\mapsto \xi_qv\in \Q$ is $C^1$.

Several Riemannian metrics on such shape spaces have been introduced and studied in the literature (see \cite{younes1998computable,ksmj04,yezzi2004metrics,mm06a,mm07,ymsm07,bauer2014constructing}), most of the times defined via reparametrization-independent positive self-adjoint differential operators, $L_q$, leading to norms taking the form
$
\normi{h}^2_{q} = \int_M h(x)^T(L_qh)(x) |q'(x)|\,dx.
$ 
For example, one can use  
\begin{equation}
\label{eq:norm.curves.0}
 \normi{h}^2_{q} = \int_M \left(\alpha\, |h(x)|^2 |q'(x)| + |h'(x)|^2|q'(x)|^{-1}\right) \,dx, 
\end{equation}
%(equivalent to the $H^1$ norm with respect to the arc-length parametrization of $q$)
 for which $L_q h = \alpha h - |q'|^{-1}\partial_x  (h'|q'|^{-1}) = \alpha h - \partial^2_s h$. Here and in the following, we use either $\partial_x f$ or $f'$ to denote the derivative of a function with respect to $x$ (assuming no ambiguity in the second case) and $\partial_s$, the arc-length derivative, denotes the operator $|q'|^{-1} \partial_x$. Higher-order derivatives have also been studied, together with the introduction of weights relying on geometric properties like length or curvature.
We will always assume that there exists a function $C$ defined on $\M$ such that $\normi{h}_q \leq C(q) \|h\|_\Q$ for all $h\in \Q$ and $q\in \M$ and that $q\mapsto \normi{h}^2_q$ is $C^1$ from $\M$ to $\mR$. 

Choosing one of these metrics, $h \mapsto \normi{h}_q$, applied to vector fields along $q$,
 we introduce the norm, applied to vector fields $v$ defined over $\mR^d$:
\begin{equation}
\label{eq:norm.curves.1}
\|v\|_{q}^2 = \lambda \|v\|_V^2 + \normi{v\circ q}^2_q.
\end{equation}
 It is important to notice that $\normi{\cdot}_q$ can be a semi-norm in this expression (i.e., one can have $\normi{h}|_q=0$ and $h\neq 0$), while still ensuring that $\|\cdot\|_q$ is a norm.  For example, one can take $\alpha = 0$ in \eqref{eq:norm.curves.0}, which has the nice property of making this semi-norm blind to translations (for which $h$ is constant). In our experiments, we use a version of this norm which is, in addition, blind to rotations, given by (as a function of the arc length)
 \begin{equation}
 \label{eq:h1.invariant}
 \normi{h}^2_{q} = \int_0^{\ell(q)} |\partial_s h|^2 \,ds - \frac{1}{\ell(q)} \left(\int_0^{\ell(q)} \partial_s h^T N(s) \,ds\right)^2
 \end{equation}
 where $\ell(q) = \int_M |\dot q| dx$ is the length of $q$ and  $N(s)$ is the unit normal to $q$ at $q(s)$.
 Another interesting norm \cite{younes1998computable,ymsm07} is 
 \[
 \normi{h}^2_{q} = \frac{1}{\ell(q)}  \int_0^{\ell(q)} |\partial_s h|^2 \,ds
 \]
 and its corresponding rotation/scale-invariant version
 \begin{multline}
 \label{eq:h1alpha.invariant}
 \normi{h}^2_{q} = \frac{1}{\ell(q)}  \int_0^{\ell(q)} |\partial_s h|^2 \,ds \\ - \left(\frac{1}{\ell(q)}\int_0^{\ell(q)} \partial_s h^T T(s) \,ds\right)^2  - \left(\frac{1}{\ell(q)}\int_0^{\ell(q)} \partial_s h^T N(s) \,ds\right)^2
 \end{multline}
 where $T(s)$ is the unit tangent. \\
 
Given such a norm, we can consider what we refer to as an {\em hybrid LDDMM} problem minimizing
\begin{equation}
\label{eq:h.lddmm}
\int_0^1 \|v(t)\|^2_{q(t)}\,dt + D(q(1), q_1) 
\end{equation}
subject to $q(0) = q_0$ and $\partial_t q(t) = v(t)\circ q(t)$.\\

There are two ways to interpret this approach. The first one, following our presentation,  is that  it modifies LDDMM by taking into account geometric properties of the curve. Alternatively, one can interpret this norm as a modification of one of the norms used on spaces of immersed curves, who generally do not prevent self intersections along geodesic paths.  From this viewpoint, the first term ($\|v\|_V^2$) is a global control ensuring the existence of a diffeomorphism of $\mR^d$ transforming the curve,  therefore guaranteeing that the curve remains embedded along any finite-energy trajectory. \\

One of the most interesting applications of this formulation is that it is easy to  generalize it to the comparison of multiple curves, say $q_1, \ldots, q_n$, using
\begin{equation}
\label{eq:multi.curves}
\|v\|_{q}^2 = \lambda \|v\|_V^2 + \normi{v\circ q_1}^2_{q_1} + \cdots + \normi{v\circ q_n}^2_{q_n}.
\end{equation}

Assuming that the curves do not overlap to start with, they will remain apart over any trajectory. However, if one chooses a small scale for the kernel of $V$, these curves will have almost no interaction unless they come close to each other. The terms $\normi{\cdot}_{q_i}$ control the shape variations of each curve separately, while $\|v\|_V$ ensures global consistency via a diffeomorphic transformation. With such a model, it becomes possible, by playing with the permissiveness of the $V$-norm relative to the curve metrics,  to transform sets of curves while ensuring that each curves evolves in an almost rigid way while their relative position in space may vary greatly. This is illustrated by a several examples in section \ref{sec:exp.1}.

Note that multiple shape comparison using the LDDMM approach was also the subject of \cite{arguillere2016registration}. The approach in that work differs from what we are proposing here, because in \cite{arguillere2016registration}, each curve was attributed its own vector field, say $v_i$, separately controlled by an RKHS norm, and global consistency was ensured by an additional vector field (similar to the $v$ that we are using here), and by equality constraints for the curve displacement, ensuring that $v\circ q = v_i\circ q$. The approach in the present paper is significantly easier to implement, avoiding, in particular, the use of constrained optimization methods.

\subsection{Maximum Principle and Optimization Algorithm}

Consider the Hamiltonian defined for $(p,q) \in \Q^*\times \Q$ and $v\in V$, by
\[
H_v(p, q) = \lform{p}{v\circ q} - \frac12 \|v\|_q^2.
\]
If we add the assumption that $q\mapsto D(q, q_1)$ is differentiable from $\M$ to $\mR$, one can prove that the 
Pontryagin's maximum principle (PMP) applies. This principle states that solutions of \eqref{eq:h.lddmm} are such that there exists $p: [0,1]\to \Q^*$ satisfying
\begin{equation}
\label{eq:pmp}
\left\{
\begin{aligned}
\partial_t q &= \partial_p H_{v(t)}(p,q)\\
\partial_t p &= -\partial_q H_{v(t)}(p,q)\\
v(t) &= \mathrm{argmax}_{w} H_w(p(t),q(t))
\end{aligned}
\right.
\end{equation}
with the boundary conditions $q(0) = q_0$ and $p(1) = - \partial_q D(q,q_1)$.  The validity of of the principle can be derived from the differentiability of the equation $\partial_t q = v\circ q$ with respect to the control ($v$) and several applications of the chain rule. (We skip the proof.)

When $\|v\|_q$ takes the form in \eqref{eq:norm.curves.1}, the Hamiltonian, considered as a function of $v$ takes the form $\lambda \|v\|_V^2 + F_{p,q}(\xi_q v)$ for a $C^1$ function $F_{p,q}: \Q \to \mR$. As a consequence, the optimal $v$ in the third equation of \eqref{eq:pmp} is such that $\lambda \scp{v}{h}_V + \lform{dF_{p,q}}{\xi_q h} = 0$ for all $h\in V$, and writing
\[
\lform{dF_{p,q}}{\xi_q h} = \lform{\xi_q^*dF_{p,q}}{h} = \scp{K\xi_q^*dF_{p,q}}{h}_V
\]
we see that $v$ should take the form $v = K\xi_q^* \alpha$ for some $\alpha\in \Q^*$. One can therefore apply the same reduction as the one which is typically used with standard LDDMM \cite{arguillere2014shape}, using $\alpha$ as a new control and minimizing
\begin{equation}
\label{eq:h.lddmm.rdc}
\frac12 \int_0^1 \|\alpha(t)\|^2_{q(t)}\,dt + D(q(1), q_1) 
\end{equation}
subject to $q(0) = q_0$ and $\partial_t q(t) = \K_{q(t)}\alpha(t)$, with $\K_q = \xi_qK\xi_q^*$ and 
\[
\|\alpha\|^2_{q} = \lambda \lform{\alpha}{\K_q\alpha} + \normi{\K_q\alpha}_q^2.
\]
The PMP can then be rewritten starting from the Hamiltonian
\[
H_\alpha(p,q) = \lform{p}{\K_q\alpha} - \frac12 \|\alpha\|_q^2.
\]
In the case we are considering in this paragraph, for which $\xi_qv = v\circ q$, $\K_q$ is given by
\[
(\K_q\alpha)(x) = \lform{\alpha}{K(q(x), q(\cdot))}.
\]
\\

In our numerical implementations, in which $q$ is represented as polygon with vertexes $(x_1, \ldots, x_N)$,  $\normi{h}_q^2$ is approximated using finite differences, so that the approximation is a function of $(h(x_i), i=1, \ldots N)$ and  of $(x_1, \ldots, x_N)$. As a result, the optimal control problem is reduced to a problem where state and controls are in $(\mR^d)^N$ (and the metric is Riemannian). The numerical results that we provide are based on this approximation (and a time approximation using a standard Euler scheme). We also recall that the computation of the  gradient of the objective function (considered as a function of the control) can be based on the PMP, using the adjoint algorithm that first computes $q(\cdot)$ by solving the first equation in \eqref{eq:pmp} starting with $q_0$ and using the control at which the gradient is computed; then sets $p(1) = - \partial_q D(q,q_1)$ before solving the second equation backward in time to obtain $p(\cdot)$; and finally computes the differential of the objective function, which is given by $\partial_v H_v(p,q)$ (or $\partial_\alpha H_\alpha(p,q)$). The gradient itself is defined by $\K_q^{-1}\partial_\alpha H_\alpha(p,q)$, which requires no operator inversion, because
\[
\K_q^{-1} \partial_\alpha H_\alpha(p,q) = p - \alpha - \G_q\K_q\alpha
\]
where $\G_q: \Q \to \Q^*$ is the operator defined by $\scp{h_1}{h_2}_q = \lform{\G_qh_1}{h_2}$.

\subsection{Experiments}
\label{sec:exp.1}
\subsubsection*{Cost function}
The end-point cost we used for our experiments is a version of the varifold metric introduced in \cite{charon2013varifold}. More precisely, we let
\begin{equation}
\label{eq:d.varifold}
D(q,q_1) = \|q\|_\chi^2 - 2\scp{q}{q_1}_\chi + \|q_1\|^2_\chi
\end{equation}
where
\begin{equation}
\label{eq:scp.varifold}
\scp{q}{q_1}_\chi^2 = \int_M\int_M \chi(q(u), q_1(u_1)) (1 + c (\nu(u)^T\nu_1(u_1))^2) |q'(u)|\,|q'_1(u_1)|\, du\,du_1,
\end{equation}
where $\nu$ and $\nu_1$ denote the unit normals to $q$ and $q_1$, $\chi$ is a Gaussian kernel
\[
\chi(x,y) = \exp(-|x-y|^2/2\tau^2),
\]
$\tau$ and $c$ being fixed parameters (we used $\tau=2$ and $c=1$ in our experiments). Because this cost function is bi-invariant by reparametrization ($D(q, q\circ \beta) = 0$ if $\beta$ is a diffeomorphism of $M$) and $\normi{\cdot}_q$ is invariant too ($\normi{h\circ \beta}_{q\circ\beta} = \normi{h}_q$), the problem is reparametrization-invariant (replacing $q_1$ by $q_1\circ \beta$ does not change the solution). 
%In all our experiments, $\normi{\cdot}_q$ was chosen as in \eqref{eq:h1alpha.invariant}.

\subsubsection{Smoothed Cardioids}
\begin{figure}
\centering
\begin{subfigure}[t]{0.3\textwidth}
\centering
\includegraphics[trim=0cm 3cm 0cm 0cm,clip,width=\textwidth]{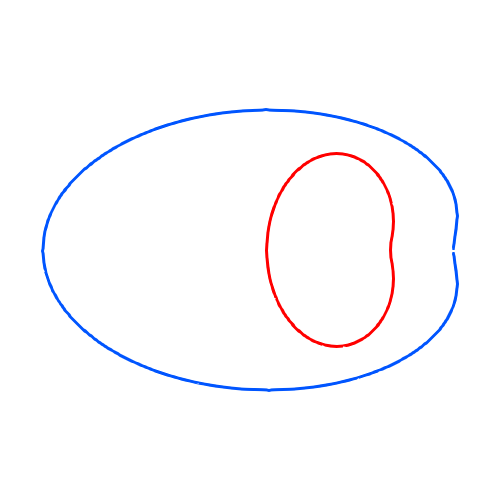}
\caption*{(a)}
\end{subfigure}
\begin{subfigure}[t]{0.3\textwidth}
\centering
\includegraphics[trim=2cm 2cm 2cm 2cm,clip,width=\textwidth]{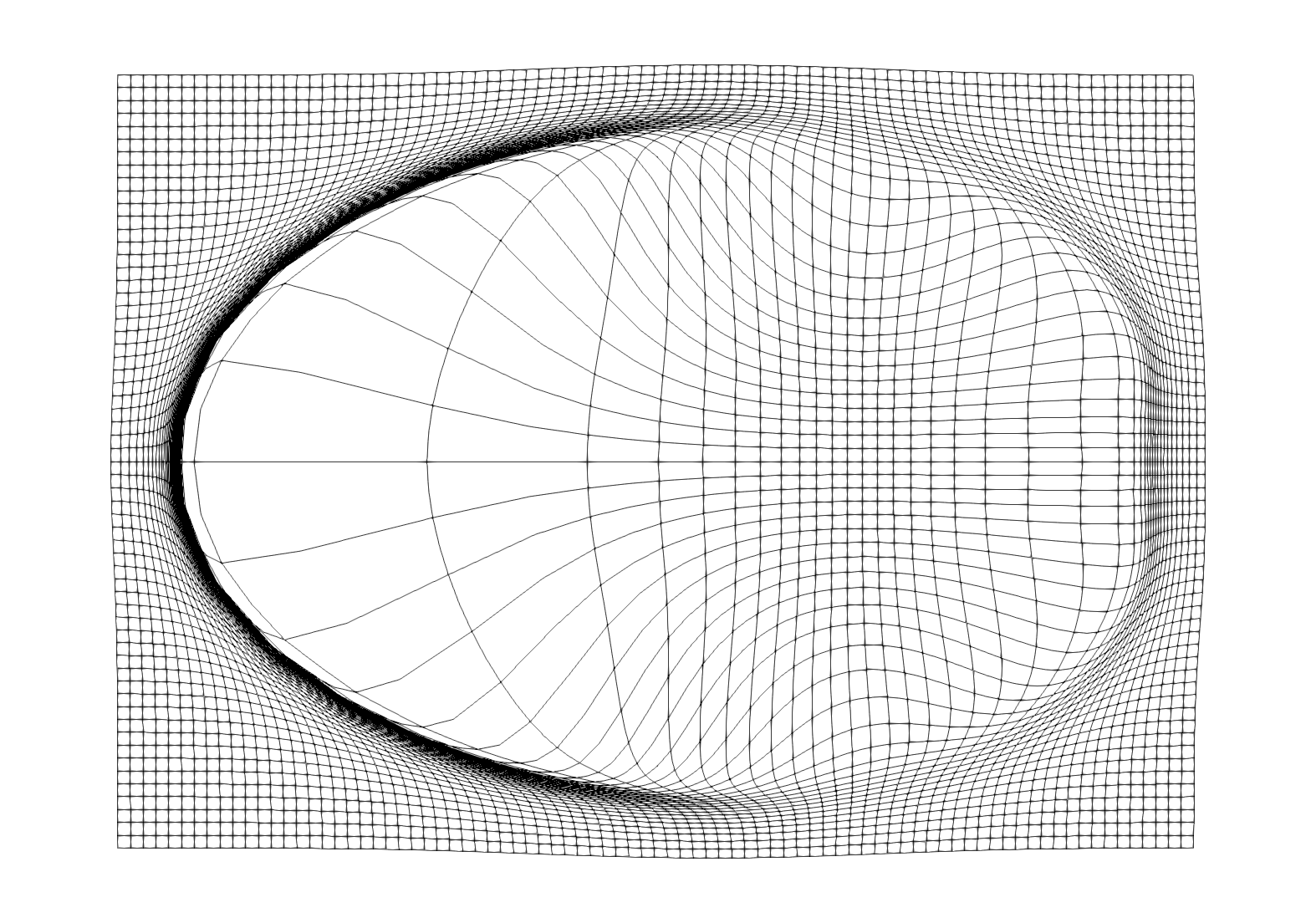}
\caption*{(b)}
\end{subfigure}
\begin{subfigure}[t]{0.3\textwidth}
\centering
\includegraphics[trim=1cm 1cm 1cm 1cm,clip,width=\textwidth]{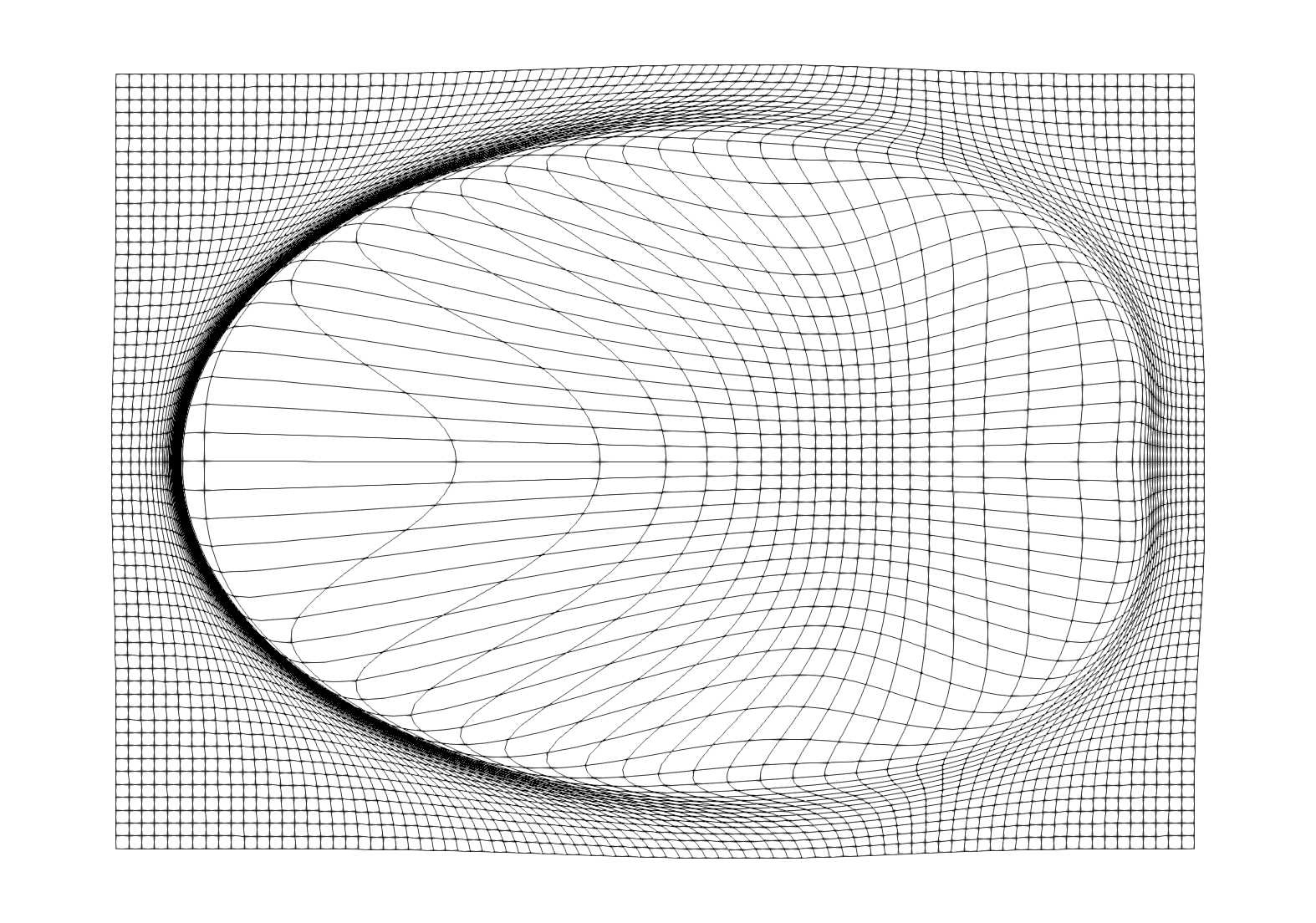}
\caption*{(c)}
\end{subfigure}
\caption{\label{fig:card.def} Smoothed Cardioids: Estimated deformations. (a) Template (red) and target (blue); (b) Standard LDDMM; (c) Hybrid LDDMM.}
\end{figure}

\begin{figure}
\centering
\begin{subfigure}[t]{0.22\textwidth}
\centering
\includegraphics[trim=1cm 1cm 1cm 1cm,clip,width=\textwidth]{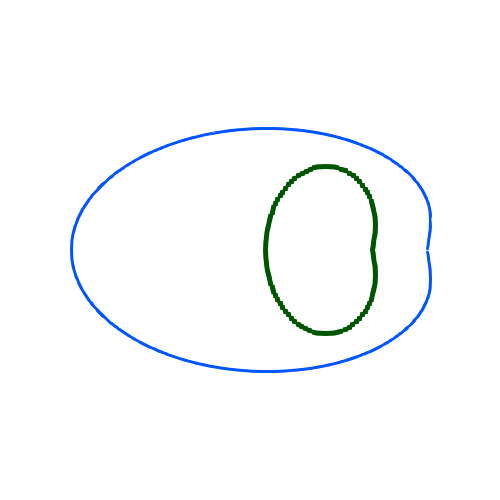}
\caption*{$t=0$}
\end{subfigure}
\begin{subfigure}[t]{0.22\textwidth}
\centering
\includegraphics[width=\textwidth]{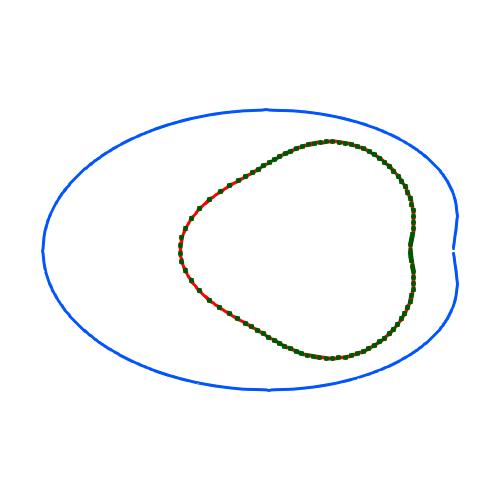}
\caption*{$t=0.3$}
\end{subfigure}
\begin{subfigure}[t]{0.22\textwidth}
\centering
\includegraphics[width=\textwidth]{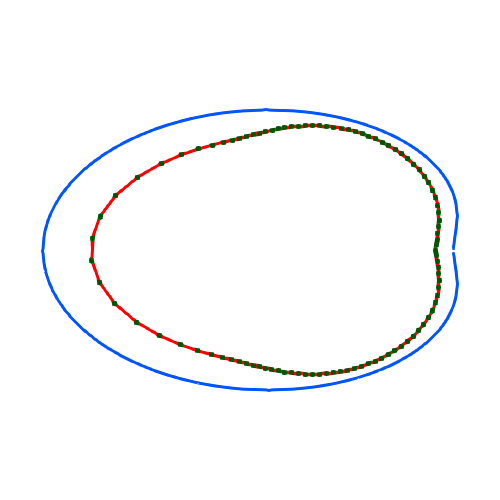}
\caption*{$t=0.7$}
\end{subfigure}
\begin{subfigure}[t]{0.22\textwidth}
\centering
\includegraphics[width=\textwidth]{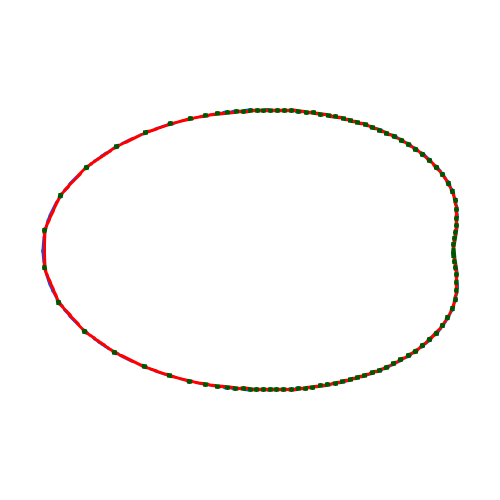}\\
\caption*{$t=1.0$}
\end{subfigure}
\begin{subfigure}[t]{0.22\textwidth}
\centering
\includegraphics[trim=1cm 1cm 1cm 1cm,clip,width=\textwidth]{cardioidst0.png}
\caption*{$t=0$}
\end{subfigure}
\begin{subfigure}[t]{0.22\textwidth}
\centering
\includegraphics[trim=1cm 1cm 1cm 1cm,clip,width=\textwidth]{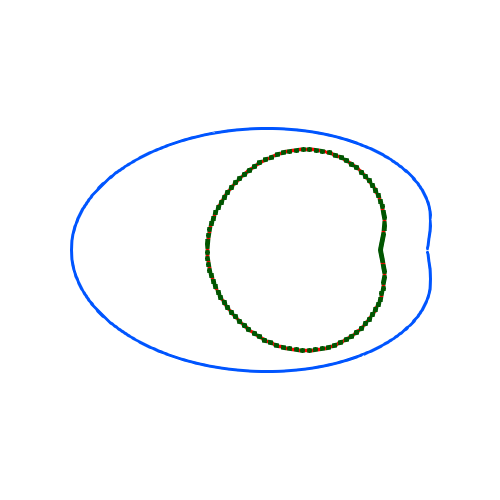}
\caption*{$t=0.3$}
\end{subfigure}
\begin{subfigure}[t]{0.22\textwidth}
\centering
\includegraphics[trim=1cm 1cm 1cm 1cm,clip,width=\textwidth]{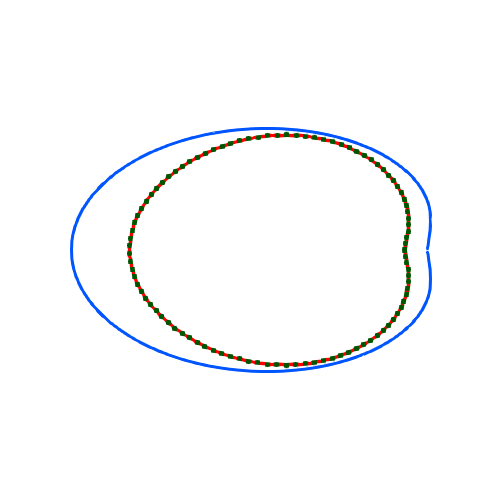}
\caption*{$t=0.7$}
\end{subfigure}
\begin{subfigure}[t]{0.22\textwidth}
\centering
\includegraphics[trim=1cm 1cm 1cm 1cm,clip,width=\textwidth]{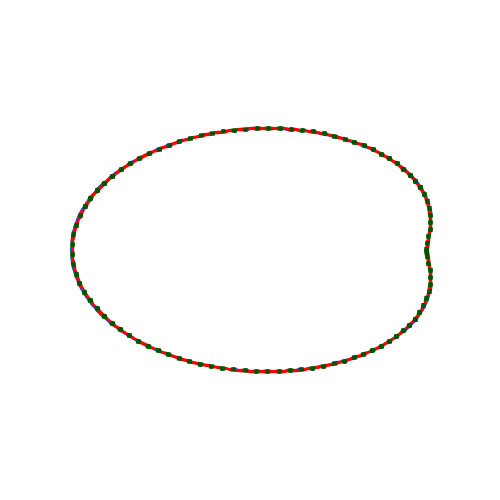}\\
\caption*{$t=1.0$}
\end{subfigure}
\caption{\label{fig:card.geod} Smoothed Cardioids: geodesics. First row: LDDMM. Second row: Hybrid LDDMM. The initial discretization of the template is uniform. It remains close to uniform along the hybrid LDDMM geodesic path, while point are spacing out in the left side and accumulating on the right side along LDDMM geodesics.}
\end{figure}

We first illustrate the impact of the additional energy term with a simple example in which two smoothed curves are registered (see right panel in Figure \ref{fig:card.def}). We used standard LDDMM with a kernel size $a=.2$ (the size of the long axis of the large cardiod being $d=10$) and hybrid LDDMM with the same kernel and $\normi{\cdot}_q$ given by \eqref{eq:h1alpha.invariant}. Both approaches perfectly align the template to the target, but their solutions differ. The LDDMM trajectories  exhibit a typical behavior in which points tend to space out during motion; this behavior is not observed in the hybrid LDDMM trajectories, because \eqref{eq:h1alpha.invariant} penalizes changes of parametrization. This can be seen in Figure \ref{fig:card.geod}, in which the deforming template is plotted in red along a geodesic path, with green dots marking the discretized points (the same color code being used in subsequent figures). The difference between the estimated registrations can also be appreciated in the last two panels of Figure \ref{fig:card.def}. Here, and in the following experiments, we used $\lambda=1$ in $\|v\|_q$, and added a multiplicative factor (between 200 and 500) in front of $\normi{\cdot}_q$ when running the hybrid version.


\subsubsection{Nested Ellipses}
\begin{figure} 
\centering
\begin{subfigure}[t]{0.3\textwidth}
\centering
\includegraphics[trim=0cm 2cm 0cm 0cm,clip,width=0.9\textwidth]{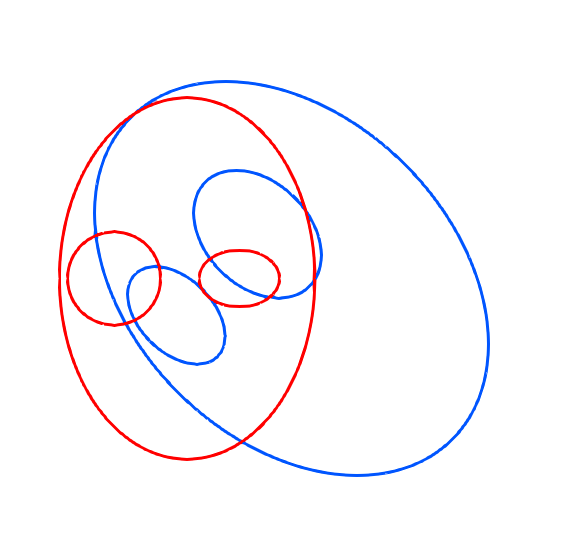}
\caption*{(a)}
\end{subfigure}
\begin{subfigure}[t]{0.3\textwidth}
\centering
\includegraphics[width=\textwidth]{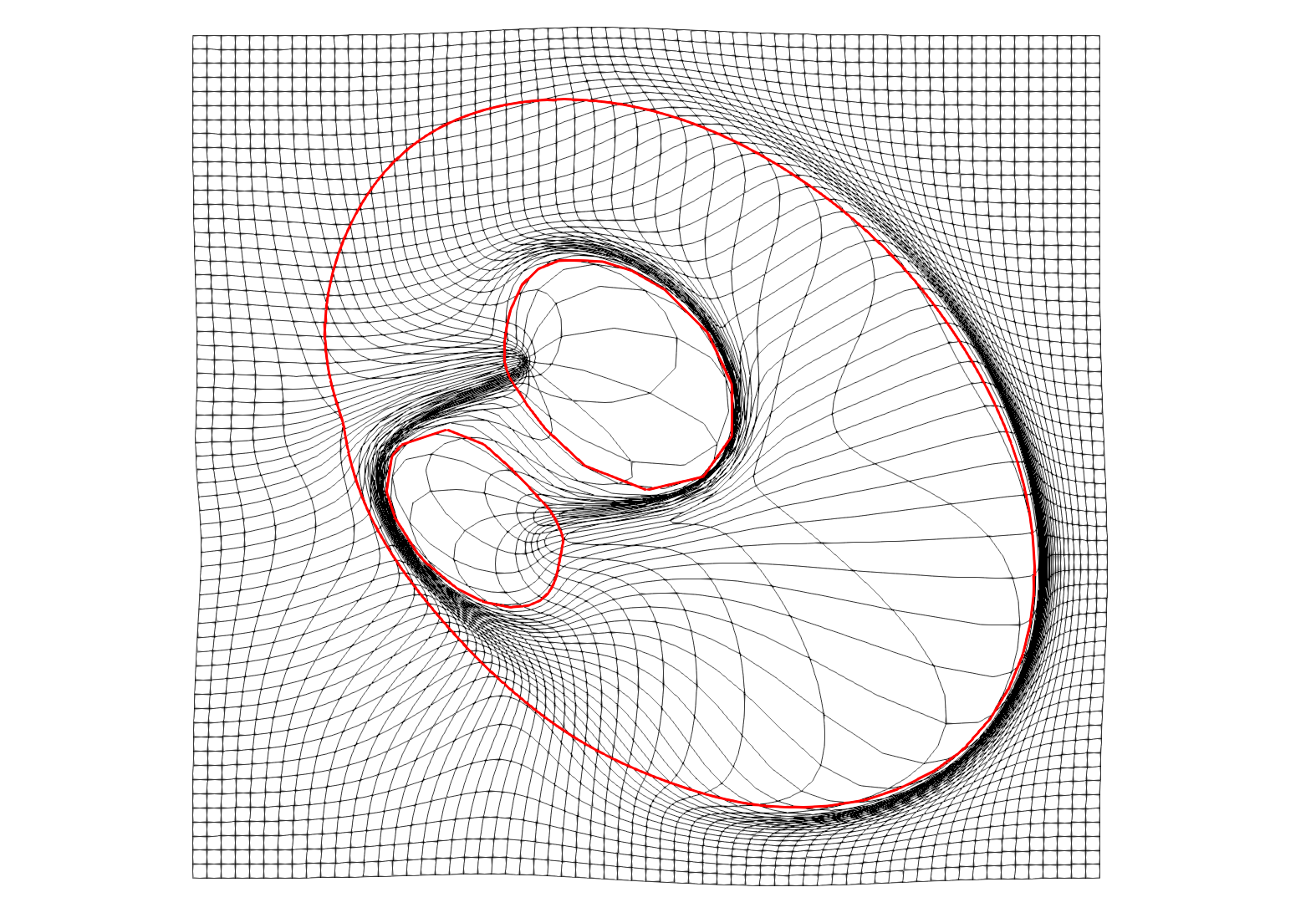}
\caption*{(b)}
\end{subfigure}\\
\begin{subfigure}[t]{0.3\textwidth}
\centering
\includegraphics[width=\textwidth]{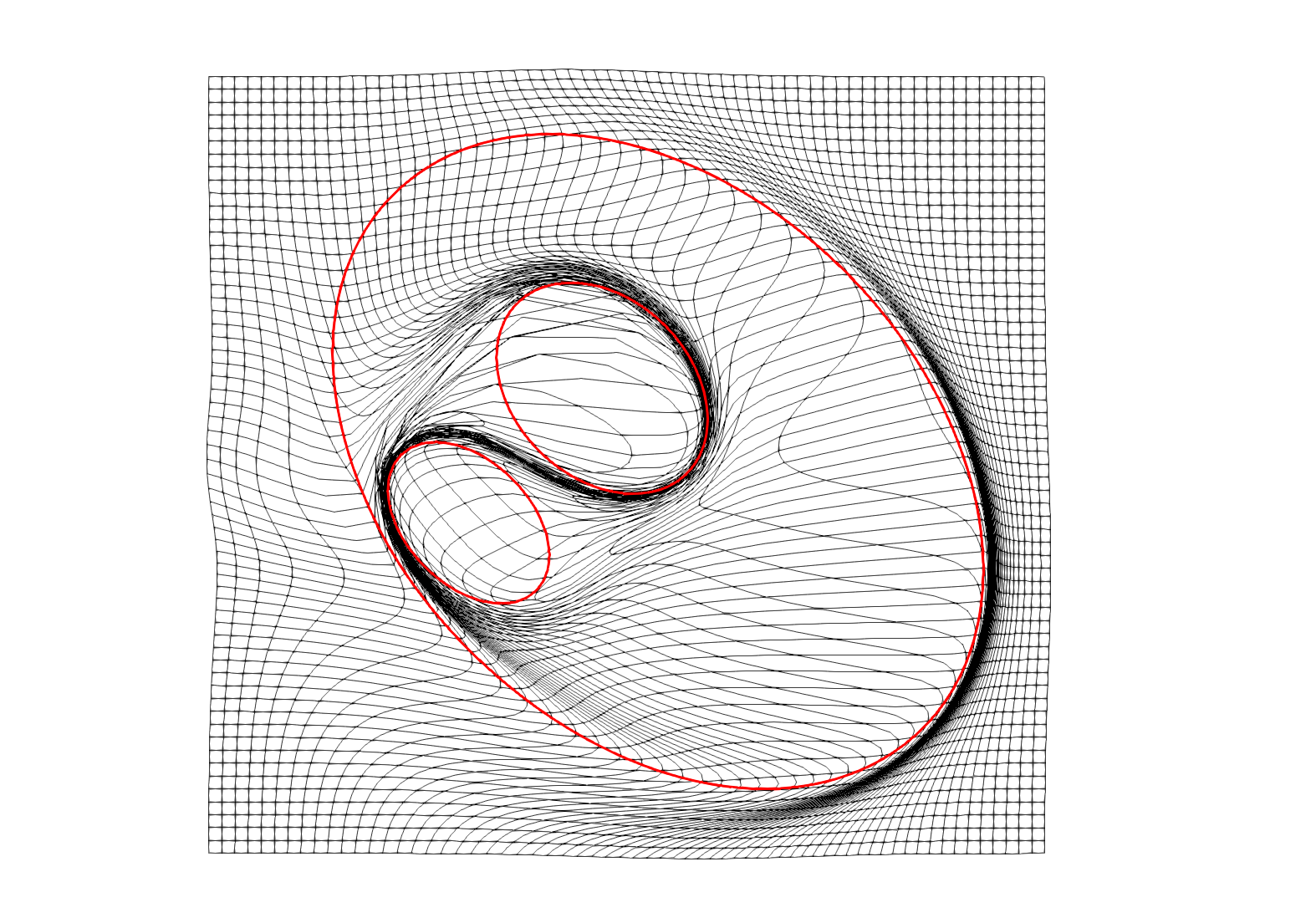}
\caption*{(c)}
\end{subfigure}
\begin{subfigure}[t]{0.3\textwidth}
\centering
\includegraphics[width=\textwidth]{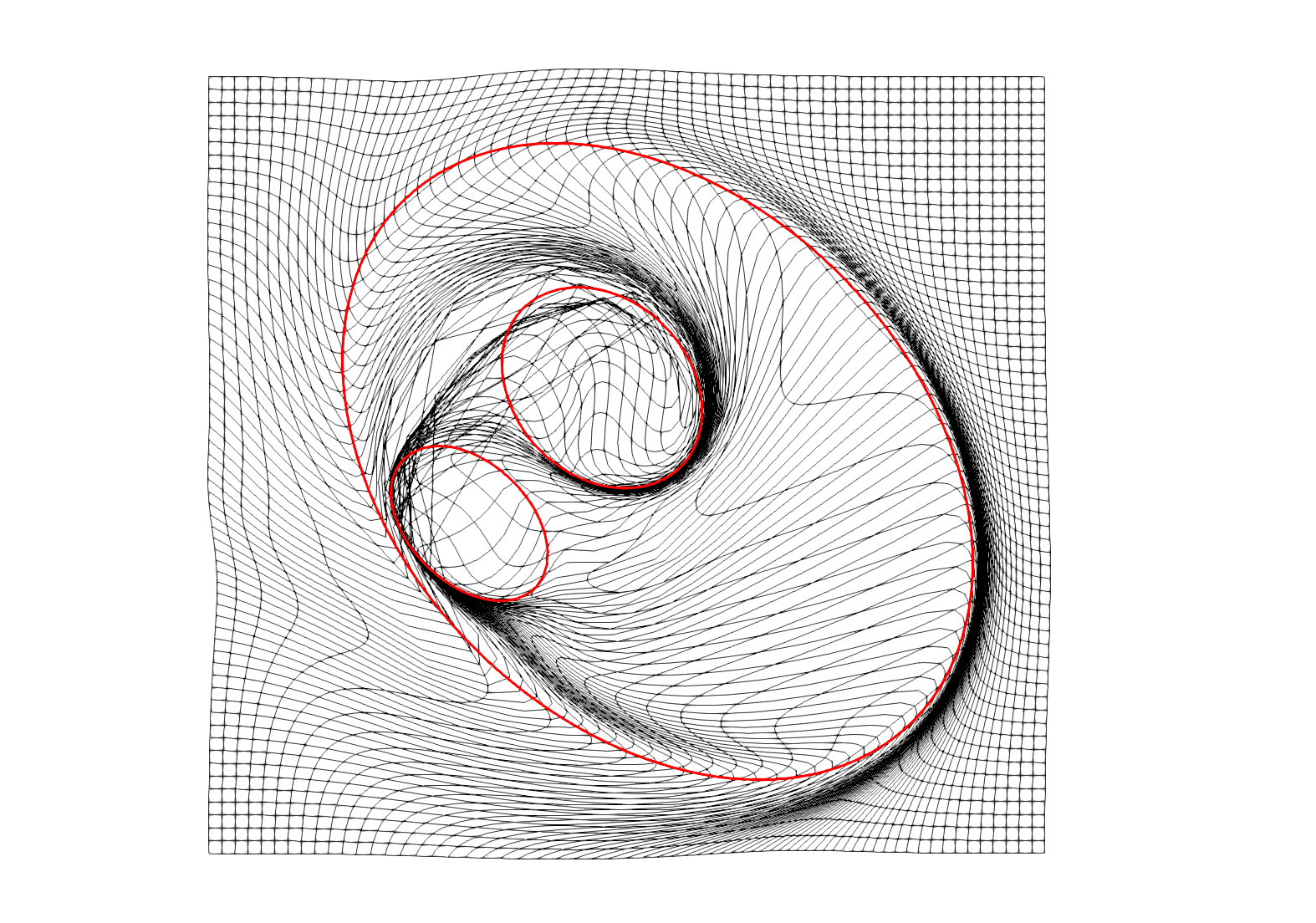}
\caption*{(d)}
\end{subfigure}
\caption{\label{fig:ellps.def} Nested Ellipses: Estimated deformations. (a) Template (red) and target (blue); (b) Standard LDDMM; (c) Hybrid LDDMM with \eqref{eq:h1alpha.invariant}. (d) Hybrid LDDMM with \eqref{eq:h1.invariant} Note that grid lines in the latter case are crossing over in this last panel. This is due to the resolution of the discretization used for this illustration compared to the size of the estimated deformation, which is diffeomorphic, as implied by the curves remaining non-intersection in the geodesic path.}
\end{figure}

Our second example is more challenging and involves multiple curves. Both template and target are composed with two small ellipses included in a large one (see Figure \ref{fig:ellps.def}). For registration, the large ellipses are paired with each other, while the small ellipses are switched, i.e., the one on the left in the template is paired with the one on the right in the target and vice versa. This is achieved by defining an end-point term as
\[
D(q^{\mathrm{large}},q_1^{\mathrm{large}})  + D(q^{\mathrm{left}},q_1^{\mathrm{right}}) + D(q^{\mathrm{right}},q_1^{\mathrm{left}})
\]
where $D$ is given by  \eqref{eq:d.varifold}.

The geodesics estimated with each method differ significantly and show interesting features. With standard LDDMM, we keep observing large reparametrization of each of the three curves, similar to what we observed in the previous example. The small ellipses avoid each other when changing places by flattening their shapes. We ran Hybrid LDDMM with $\normi{\cdot}_q$ given by \eqref{eq:h1alpha.invariant} and \eqref{eq:h1.invariant}. In both cases,  the reparametrization is uniform along each curve.  With \eqref{eq:h1alpha.invariant}, which is scale and rotation invariant, the small ellipses shrink when crossing each other, before growing back to match the target. When using \eqref{eq:h1.invariant} (which is only rotation invariant), shrinking is not free anymore, and the curves make a wide berth to avoid each other. The kernel width was the same in all three experiments, in which we took $a=0.2$.

\begin{figure}
\centering
\begin{subfigure}[t]{0.22\textwidth}
\centering
\includegraphics[trim=1cm 1cm 1cm 1cm,clip,width=\textwidth]{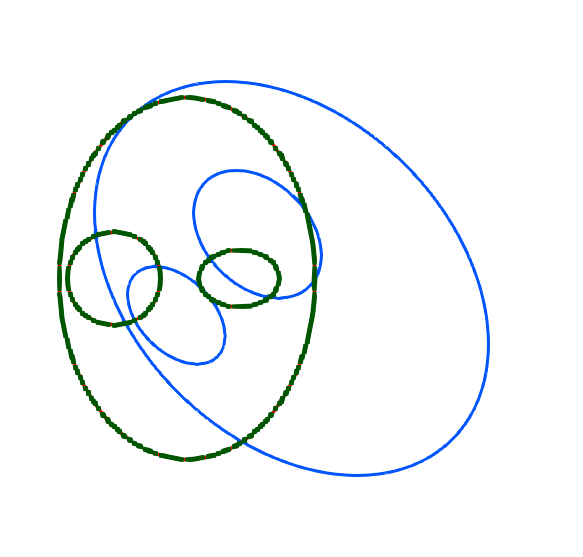}
\caption*{$t=0$}
\end{subfigure}
\begin{subfigure}[t]{0.22\textwidth}
\centering
\includegraphics[trim=1cm 1cm 1cm 1cm,clip,width=\textwidth]{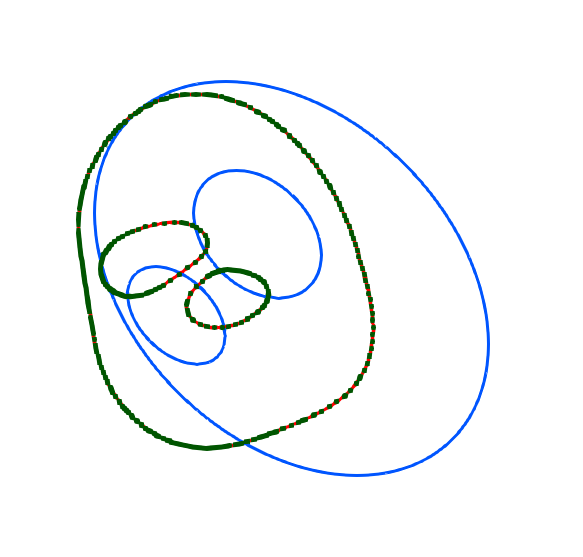}
\caption*{$t=0.3$}
\end{subfigure}
\begin{subfigure}[t]{0.22\textwidth}
\centering
\includegraphics[trim=1cm 1cm 1cm 1cm,clip,width=\textwidth]{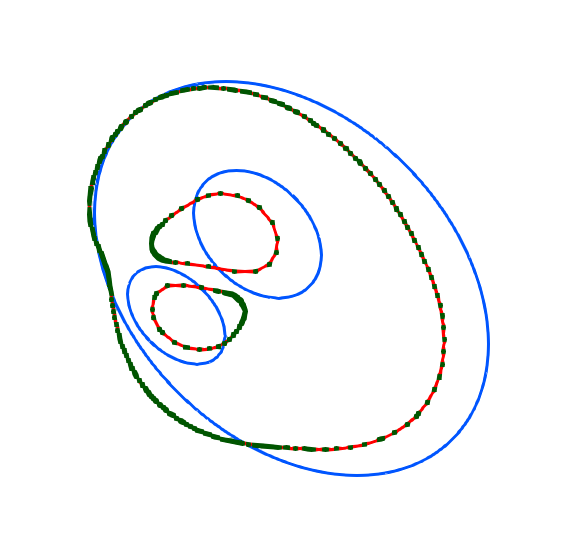}
\caption*{$t=0.7$}
\end{subfigure}
\begin{subfigure}[t]{0.22\textwidth}
\centering
\includegraphics[trim=1cm 1cm 1cm 1cm,clip,width=\textwidth]{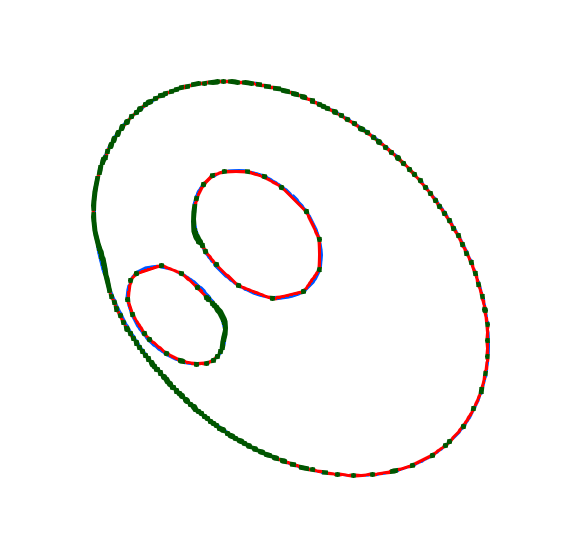}\\
\caption*{$t=1.0$}
\end{subfigure}
\begin{subfigure}[t]{0.22\textwidth}
\centering
\includegraphics[trim=1cm 1cm 1cm 1cm,clip,width=\textwidth]{ellipsest0.png}
\caption*{$t=0$}
\end{subfigure}
\begin{subfigure}[t]{0.22\textwidth}
\centering
\includegraphics[trim=1cm 1cm 1cm 1cm,clip,width=\textwidth]{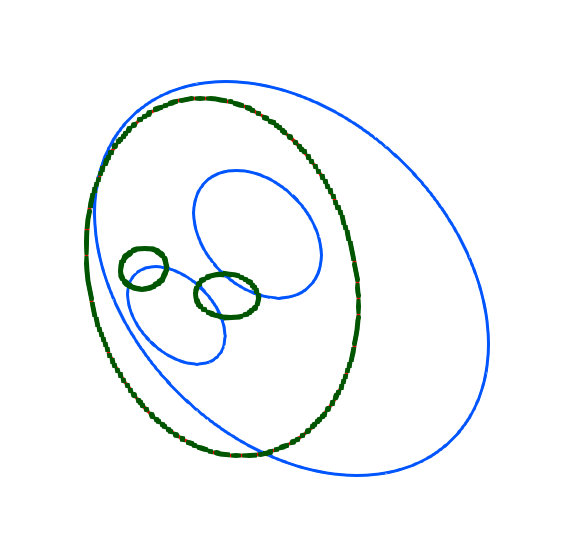}
\caption*{$t=0.3$}
\end{subfigure}
\begin{subfigure}[t]{0.22\textwidth}
\centering
\includegraphics[trim=1cm 1cm 1cm 1cm,clip,width=\textwidth]{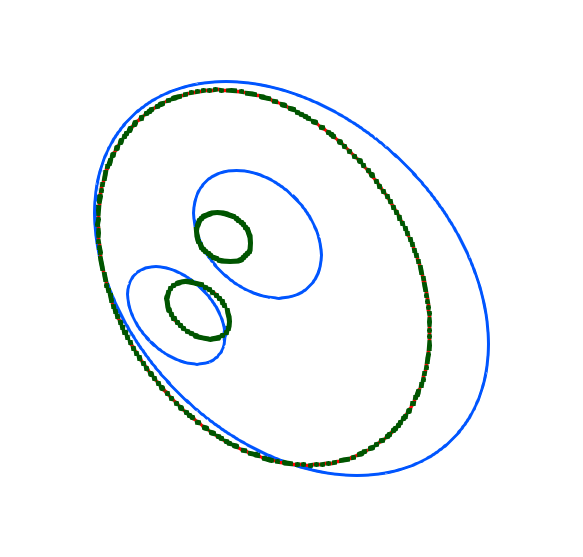}
\caption*{$t=0.7$}
\end{subfigure}
\begin{subfigure}[t]{0.22\textwidth}
\centering
\includegraphics[trim=1cm 1cm 1cm 1cm,clip,width=\textwidth]{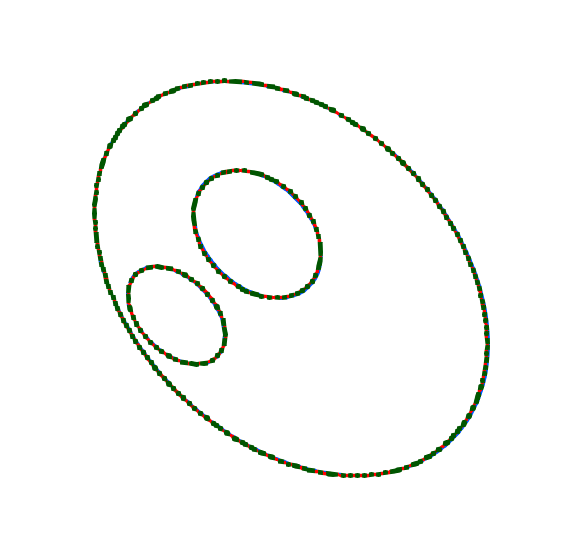}\\
\caption*{$t=1.0$}
\end{subfigure}
\begin{subfigure}[t]{0.22\textwidth}
\centering
\includegraphics[trim=1cm 1cm 1cm 1cm,clip,width=\textwidth]{ellipsest0.png}
\caption*{$t=0$}
\end{subfigure}
\begin{subfigure}[t]{0.22\textwidth}
\centering
\includegraphics[trim=1cm 1cm 1cm 1cm,clip,width=\textwidth]{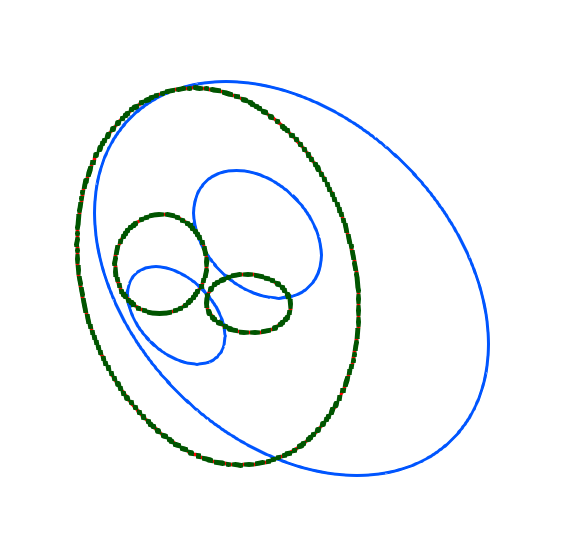}
\caption*{$t=0.3$}
\end{subfigure}
\begin{subfigure}[t]{0.22\textwidth}
\centering
\includegraphics[trim=1cm 1cm 1cm 1cm,clip,width=\textwidth]{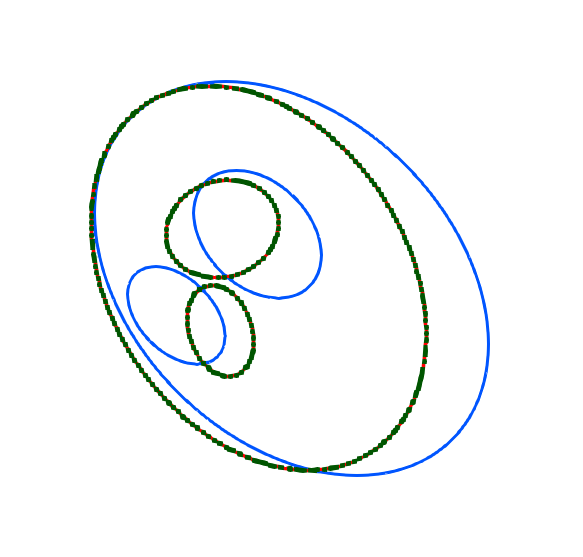}
\caption*{$t=0.7$}
\end{subfigure}
\begin{subfigure}[t]{0.22\textwidth}
\centering
\includegraphics[trim=1cm 1cm 1cm 1cm,clip,width=\textwidth]{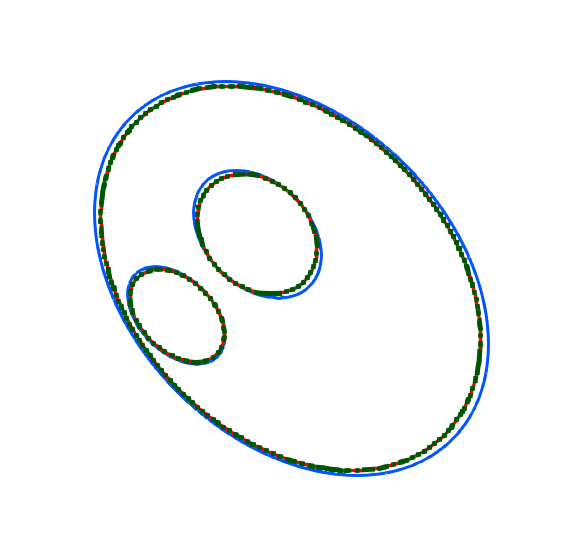}\\
\caption*{$t=1.0$}
\end{subfigure}
\caption{\label{fig:ellps.geod}Ellipses: geodesics. First row: LDDMM. Second row: Hybrid LDDMM. The initial discretization of the template is uniform.}
\end{figure}

\subsubsection{Rays}
\begin{figure} 
\centering
\begin{subfigure}[t]{0.3\textwidth}
\centering
\includegraphics[trim=0cm 4cm 0cm 0cm,clip,width=\textwidth]{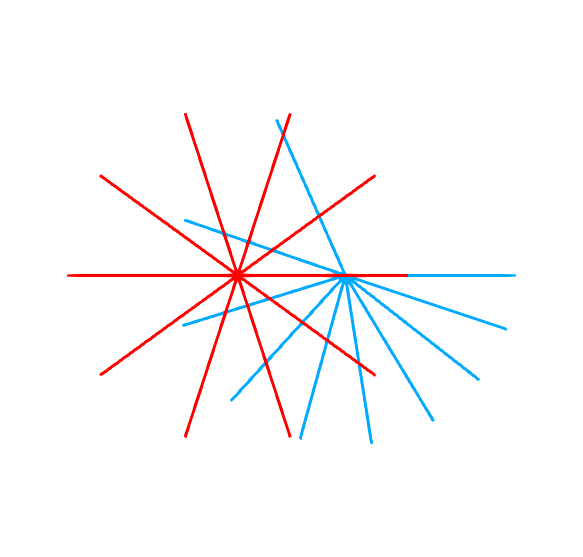}
\caption*{(a)}
\end{subfigure}
\begin{subfigure}[t]{0.3\textwidth}
\centering
\includegraphics[width=\textwidth]{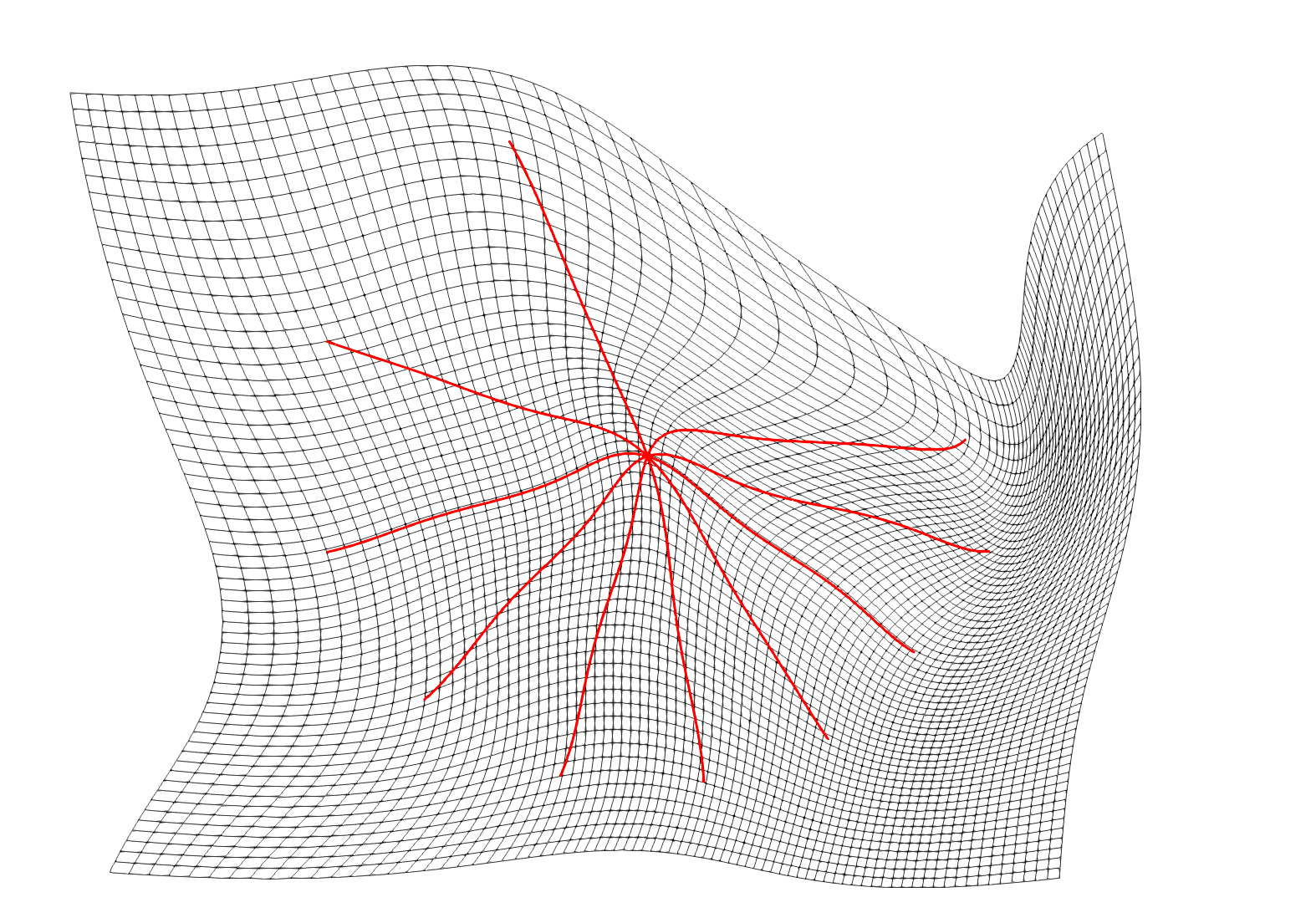}
\caption*{(b)}
\end{subfigure}
\begin{subfigure}[t]{0.3\textwidth}
\centering
\includegraphics[width=\textwidth]{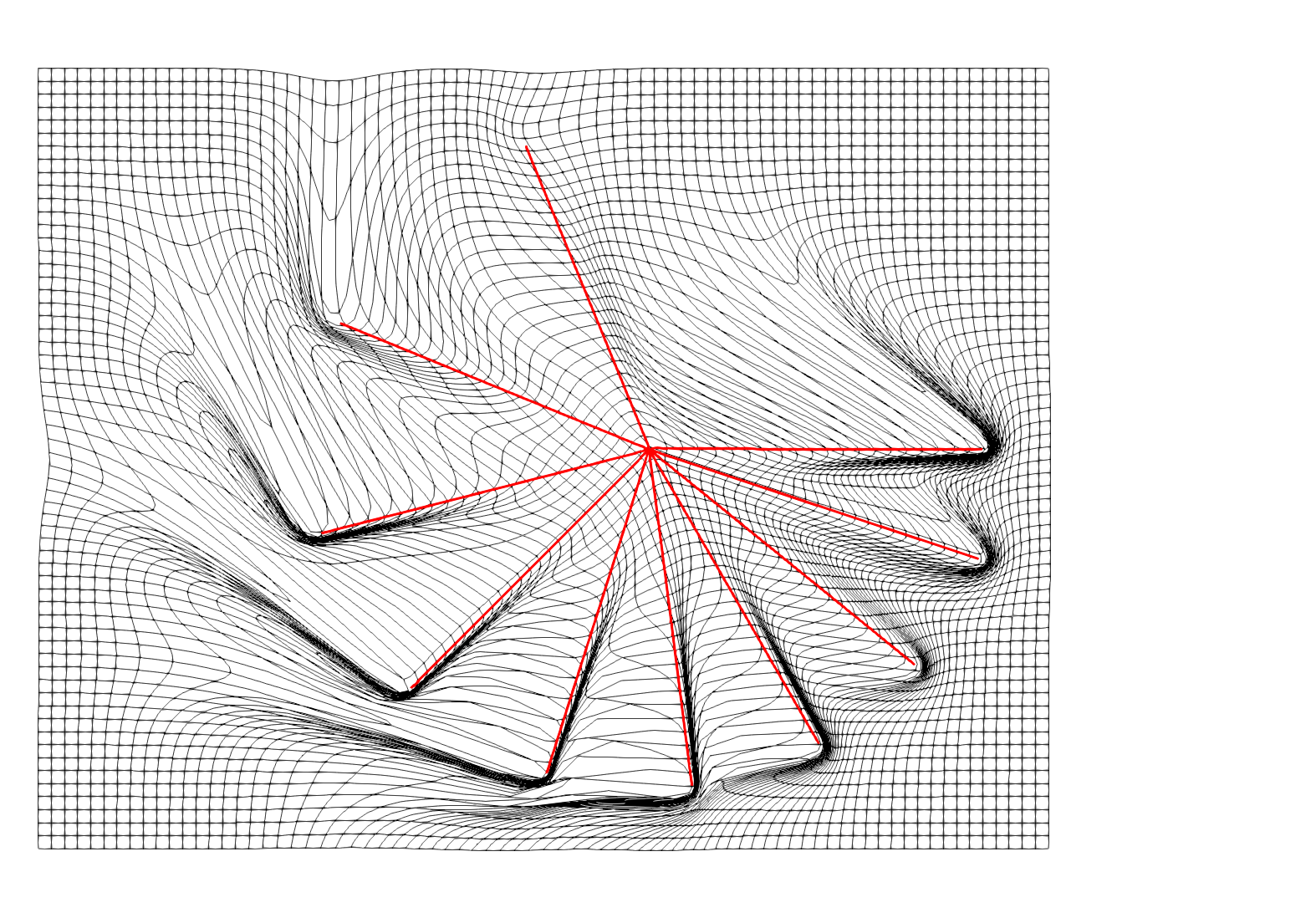}
\caption*{(c)}
\end{subfigure}
\caption{\label{fig:rays.def} Rays: Estimated deformations. (a) Template (red) and target (blue); (b) Standard LDDMM; (c) Hybrid LDDMM.}
\end{figure}

We now compare configurations of $m=10$ line segments stemming from a common origin (Figure \ref{fig:rays.def}). The segments' orientations are sampled uniformly over $[0, 2\pi]$ ($\theta_k = 2k\pi/m$, $k=0, \ldots, m-1$) in the template, but not in the target ($\theta_k  = 2\pi\sqrt{k/m}$, $k=0, \ldots, m-1$). The target is moreover slightly translated. Here, and for the examples that follow, the cost function considers the curves as unlabeled (no correspondence information is used). Formally, this corresponds to considering that the curves are parametrized over the unions of $m$ copies of $M$, $M^{(m)} = \bigcup_{k=1}^m \{k\}\times M$ and using $M^{(m)}$ in place of $M$ in \eqref{eq:scp.varifold}. This choice makes the matching problem significantly harder, creating possible local minima in the cost function. Such local minima actually trap the LDDMM algorithm when using small kernel sizes, and the solution provided in our experiments use a rather large kernel size,  $a = L/5$, where $L$ is the common lengths of the segments. The hybrid model uses $a = L/25$ combined with \eqref{eq:h1.invariant}, the $H^1$ norm  corrected for rotations. 

\begin{figure}
\centering
\begin{subfigure}[t]{0.22\textwidth}
\centering
\includegraphics[trim=1cm 1cm 1cm 1cm,clip,width=\textwidth]{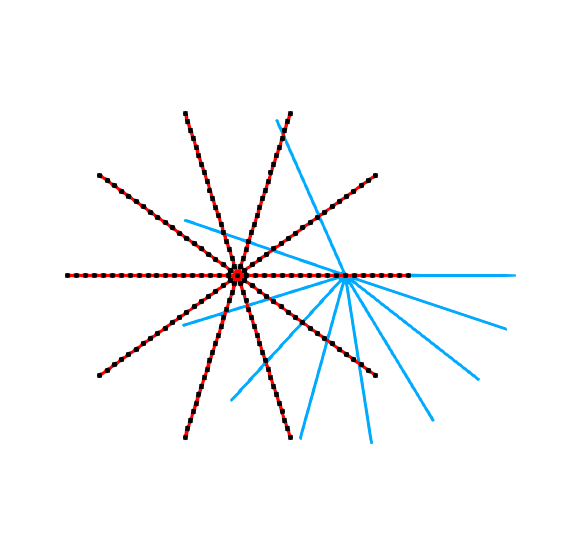}
\caption*{$t=0$}
\end{subfigure}
\begin{subfigure}[t]{0.22\textwidth}
\centering
\includegraphics[trim=1cm 1cm 1cm 1cm,clip,width=\textwidth]{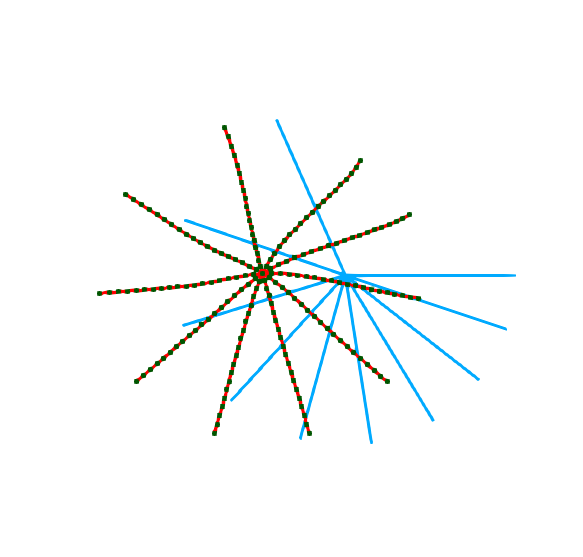}
\caption*{$t=0.3$}
\end{subfigure}
\begin{subfigure}[t]{0.22\textwidth}
\centering
\includegraphics[trim=1cm 1cm 1cm 1cm,clip,width=\textwidth]{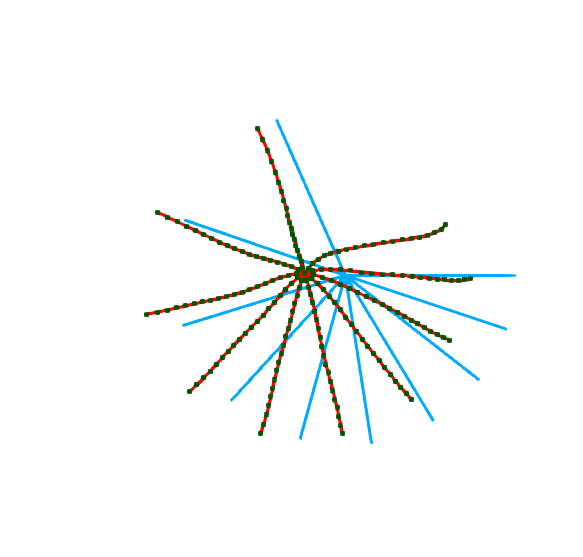}
\caption*{$t=0.7$}
\end{subfigure}
\begin{subfigure}[t]{0.22\textwidth}
\centering
\includegraphics[trim=1cm 1cm 1cm 1cm,clip,width=\textwidth]{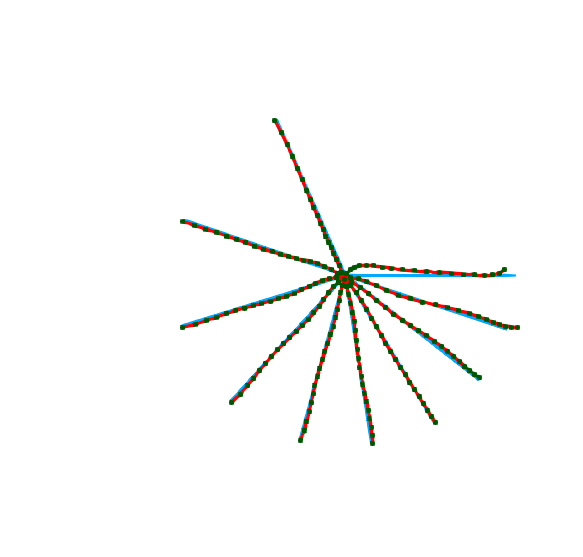}\\
\caption*{$t=1.0$}
\end{subfigure}
\begin{subfigure}[t]{0.22\textwidth}
\centering
\includegraphics[trim=1cm 1cm 1cm 1cm,clip,width=\textwidth]{rayst0.png}
\caption*{$t=0$}
\end{subfigure}
\begin{subfigure}[t]{0.22\textwidth}
\centering
\includegraphics[trim=1cm 1cm 1cm 1cm,clip,width=\textwidth]{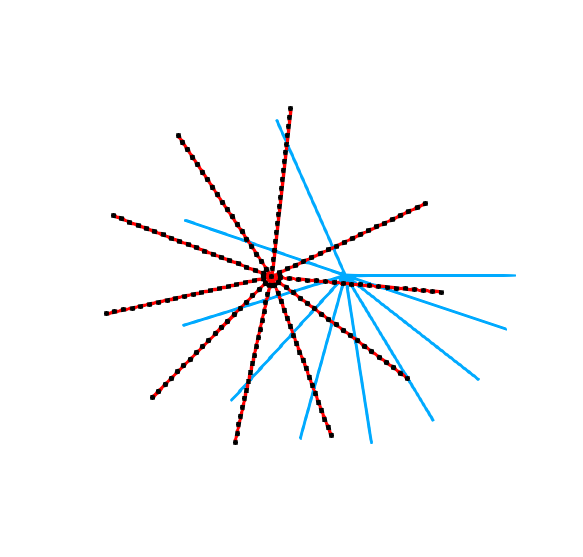}
\caption*{$t=0.3$}
\end{subfigure}
\begin{subfigure}[t]{0.22\textwidth}
\centering
\includegraphics[trim=1cm 1cm 1cm 1cm,clip,width=\textwidth]{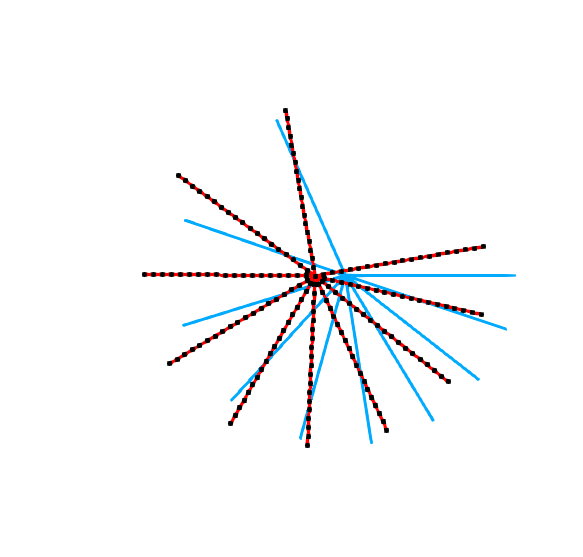}
\caption*{$t=0.7$}
\end{subfigure}
\begin{subfigure}[t]{0.22\textwidth}
\centering
\includegraphics[trim=1cm 1cm 1cm 1cm,clip,width=\textwidth]{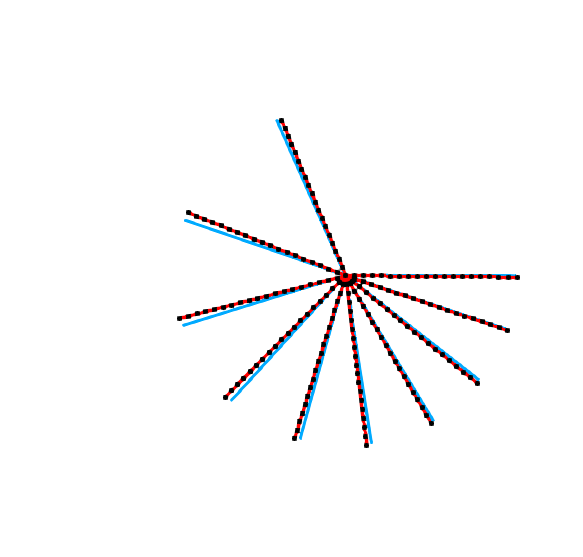}\\
\caption*{$t=1.0$}
\end{subfigure}
\caption{\label{fig:rays.geod} Rays: geodesics. First row: LDDMM. Second row: Hybrid LDDMM. The initial discretization of the template is uniform.}
\end{figure}

As a result of the use of a large kernel in the LDDMM case, the obtained solution does not achieve a perfect transformation of the first segment (the one requiring the largest rotation) which is curved at the end-point of the geodesic (see Figure \ref{fig:rays.geod}). The segments remain perfectly straight along the geodesic estimated with the hybrid norm (visually at least: an exact transformation of the rays would not be diffeomorphic, but the deviation from a straight line happens below the discretization level chosen for the curves).  The effect of the kernel size is also apparent in the estimated transformations, depicted in Figure \ref{fig:rays.def}. Similar to the previous examples, the reparametrization of the segments is more pronounced with standard LDDMM.

\subsubsection{Half Circles} 
\begin{figure} 
\centering
\begin{subfigure}[t]{0.3\textwidth}
\centering
\includegraphics[width=\textwidth]{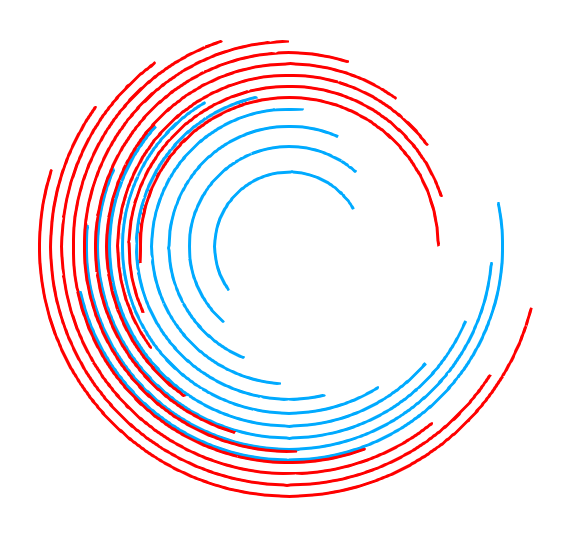}
\caption*{(a)}
\end{subfigure}
\begin{subfigure}[t]{0.3\textwidth}
\centering
\includegraphics[trim=8cm 0cm 8cm 4cm,clip,width=\textwidth]{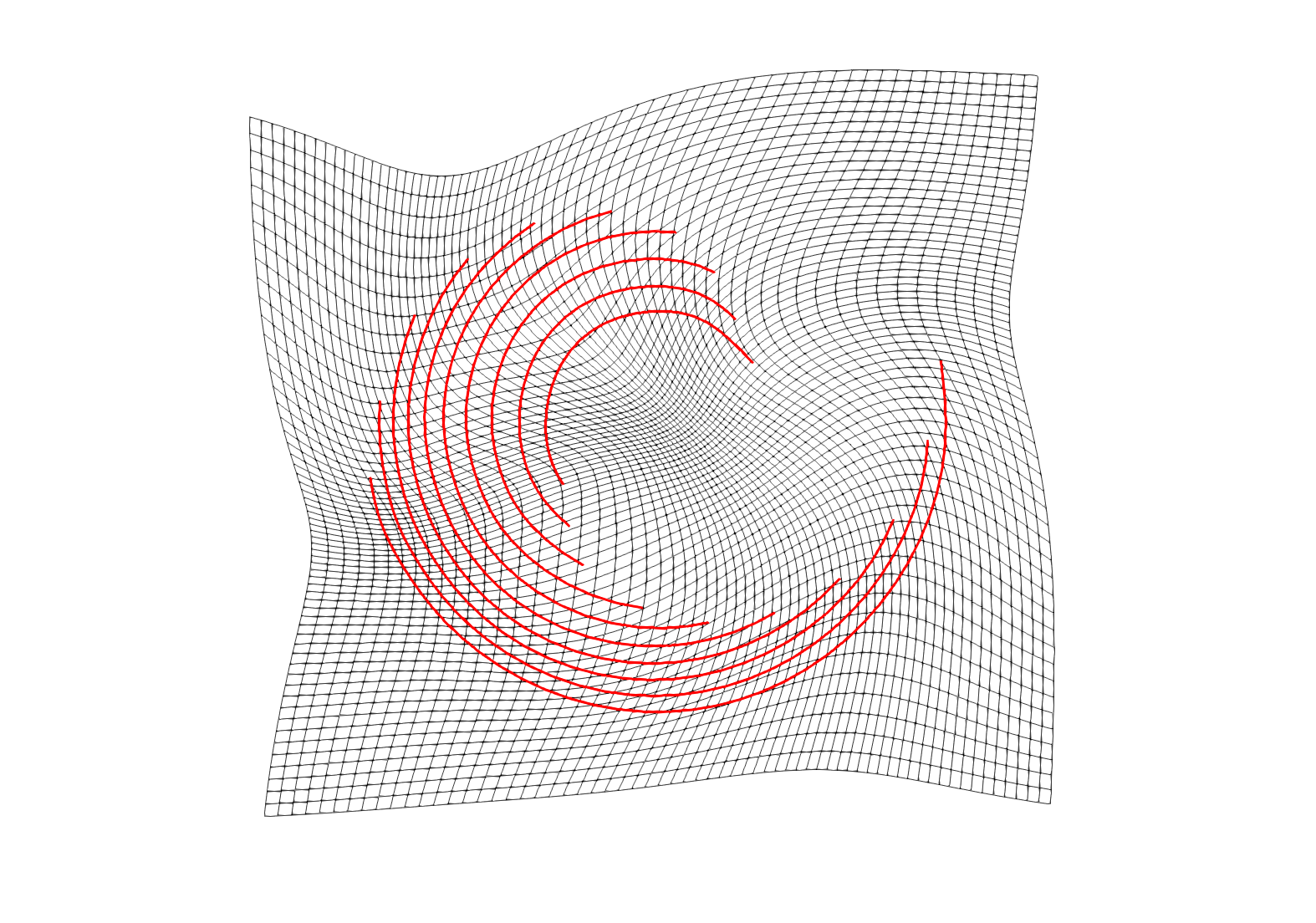}
\caption*{(b)}
\end{subfigure}
\begin{subfigure}[t]{0.3\textwidth}
\centering
\includegraphics[trim=8cm 0cm 8cm 4cm,clip,width=\textwidth]{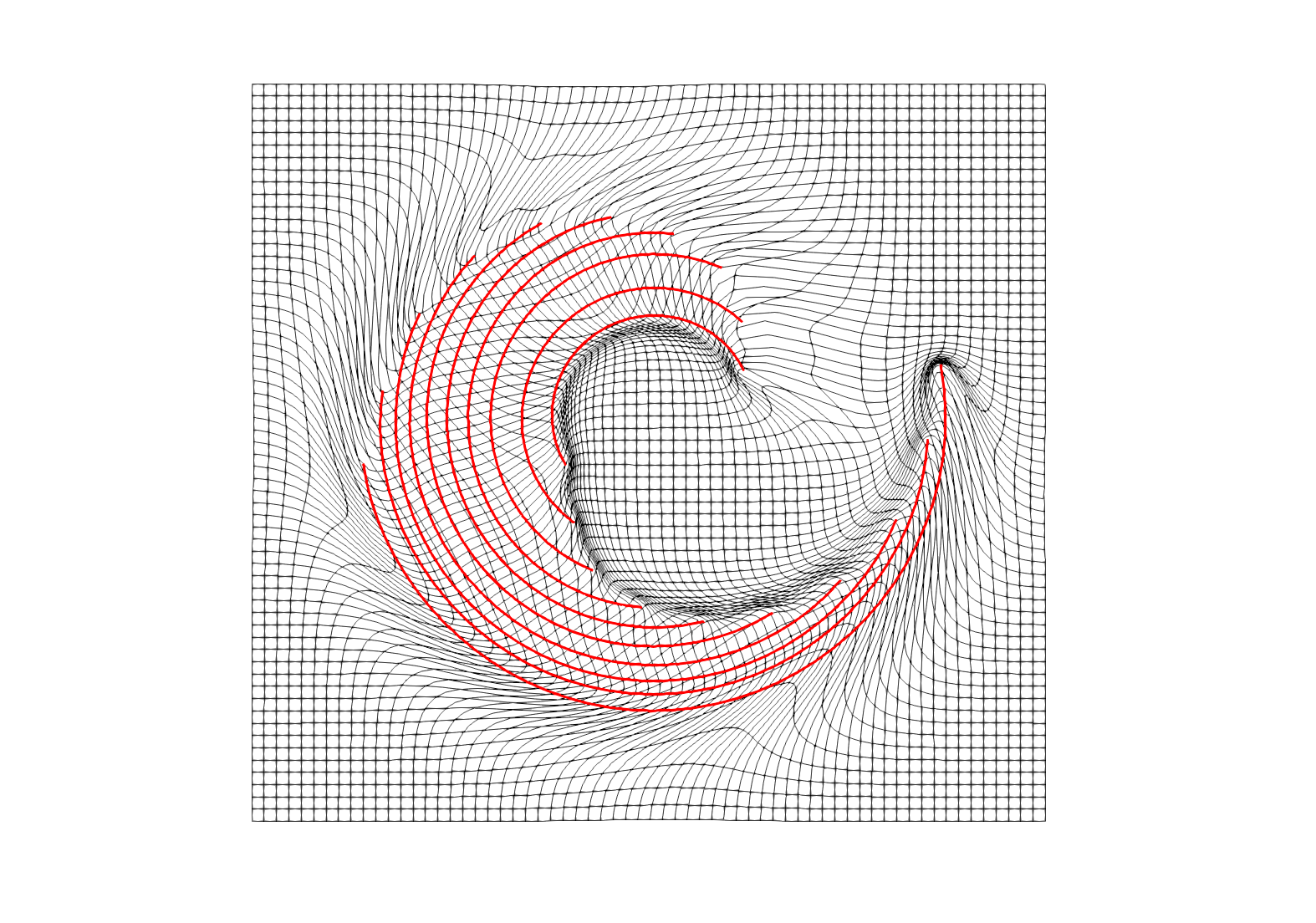}
\caption*{(c)}
\end{subfigure}
\caption{\label{fig:manyCurves.def} Arcs of circle: Estimated deformations. (a) Template (red) and target (blue); (b) Standard LDDMM; (c) Hybrid LDDMM.}
\end{figure}

Our last 2D example is similar to the previous one, using half circles, with various radii, instead of straight lines. We used the rotation- and scale-invariant norm \eqref{eq:h1alpha.invariant} in the hybrid case, and the kernel sizes were $L/5$ and $L/25$ for standard and hybrid LDDMM, $L$ being the radius of the largest circle. Both methods do a good job in registering the target to the template, but find different solutions as seen in Figure \ref{fig:manyCurves.def}. The LDDMM solution tends to compress the  space in the middle of the estimated pattern, while hybrid LDDMM estimates a motion closer to a rotation (which is cheap with the considered norm).

\begin{figure}
\centering
\begin{subfigure}[t]{0.22\textwidth}
\centering
\includegraphics[trim=1cm 1cm 1cm 1cm,clip,width=\textwidth]{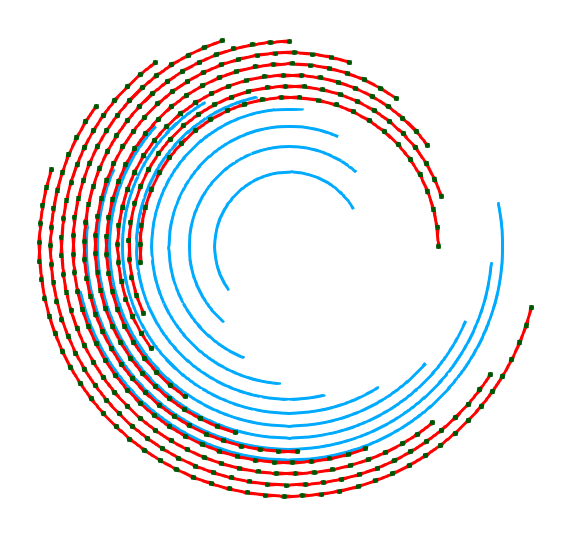}
\caption*{$t=0$}
\end{subfigure}
\begin{subfigure}[t]{0.22\textwidth}
\centering
\includegraphics[trim=1cm 1cm 1cm 1cm,clip,width=\textwidth]{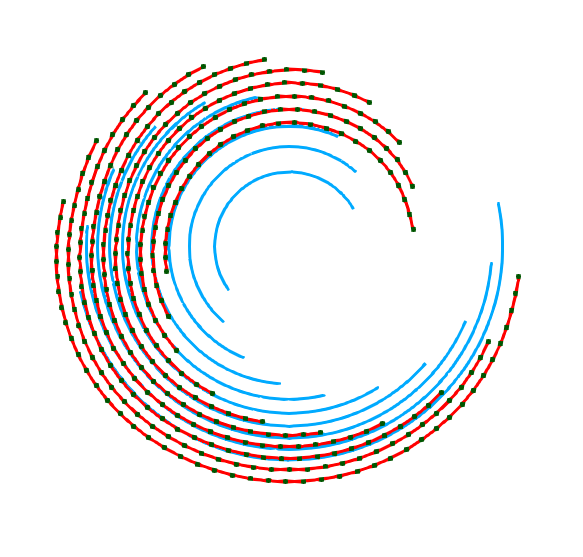}
\caption*{$t=0.3$}
\end{subfigure}
\begin{subfigure}[t]{0.22\textwidth}
\centering
\includegraphics[trim=1cm 1cm 1cm 1cm,clip,width=\textwidth]{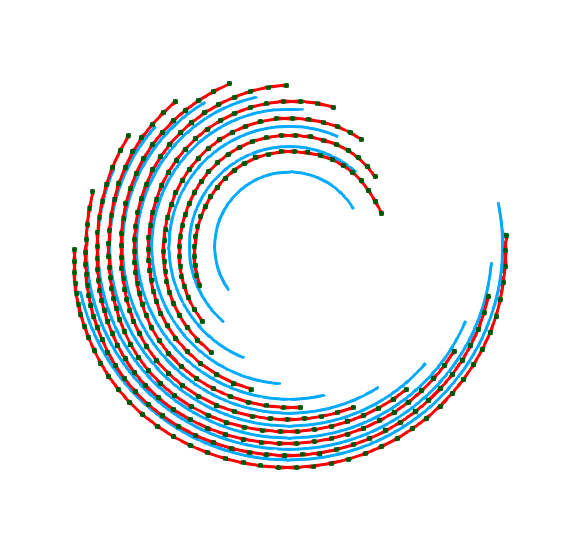}
\caption*{$t=0.7$}
\end{subfigure}
\begin{subfigure}[t]{0.22\textwidth}
\centering
\includegraphics[trim=1cm 1cm 1cm 1cm,clip,width=\textwidth]{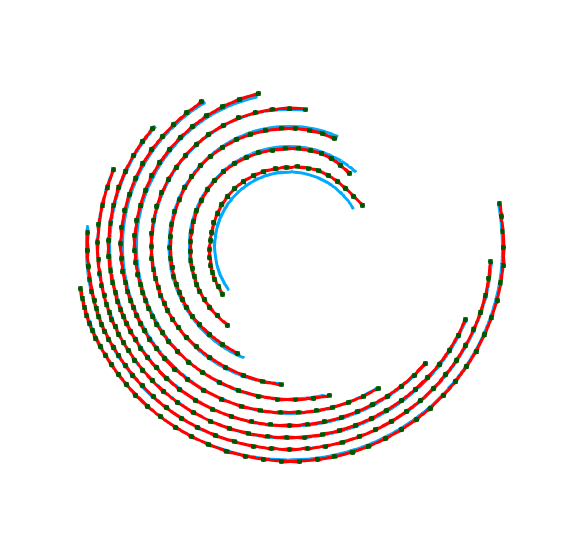}\\
\caption*{$t=1.0$}
\end{subfigure}
\begin{subfigure}[t]{0.22\textwidth}
\centering
\includegraphics[trim=1cm 1cm 1cm 1cm,clip,width=\textwidth]{manyCurvest0.png}
\caption*{$t=0$}
\end{subfigure}
\begin{subfigure}[t]{0.22\textwidth}
\centering
\includegraphics[trim=1cm 1cm 1cm 1cm,clip,width=\textwidth]{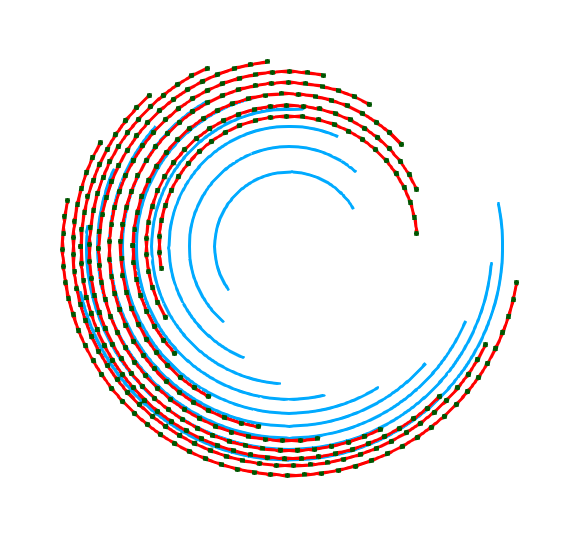}
\caption*{$t=0.3$}
\end{subfigure}
\begin{subfigure}[t]{0.22\textwidth}
\centering
\includegraphics[trim=1cm 1cm 1cm 1cm,clip,width=\textwidth]{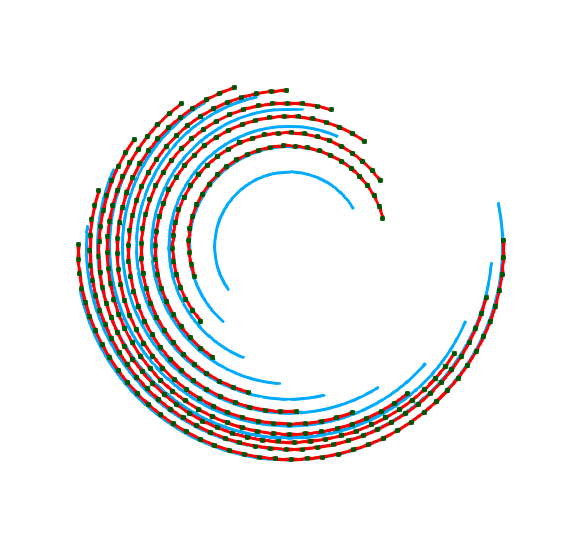}
\caption*{$t=0.7$}
\end{subfigure}
\begin{subfigure}[t]{0.22\textwidth}
\centering
\includegraphics[trim=1cm 1cm 1cm 1cm,clip,width=\textwidth]{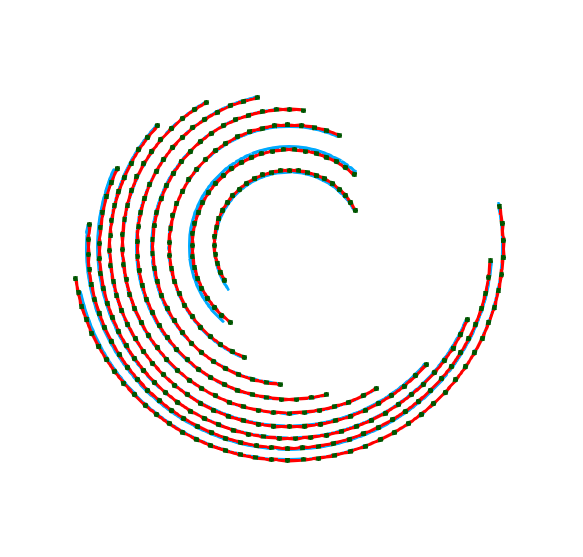}\\
\caption*{$t=1.0$}
\end{subfigure}
\caption{\label{fig:manyCurves.geod} Arcs of circle: geodesics. First row: LDDMM. Second row: Hybrid LDDMM. The initial discretization of the template is uniform.}
\end{figure}

\section{Surfaces}
The same approach can also be used with surfaces. At high level, not much needs to be modified formally from the curve case, simply letting $M$ be the unit disc, or the unit sphere, or any other  manifold sharing the topology of the considered surfaces. One can find a large collection of possible choices  for $\normi{\cdot}_q$ in \cite{bauer2011sobolev,bauer2011new,kurtek2011parameterization,bauer2014overview}. In the experiment that follows, we used the one of simplest options, letting
\begin{equation}
\label{eq:h1.surf}
\normi{h}_q^2 = \int_{S} |\nabla_S h|^2 d\sigma_S
\end{equation}
where $S = q(M)$ with area measure $\sigma_S$ and Riemannian gradient $\nabla_S$ (it is therefore important that $q$ remains an embedding at all times along finite energy paths). 

\begin{figure}
\centering
\begin{subfigure}[t]{0.4\textwidth}
\centering
\includegraphics[trim=4cm 2cm 3cm 0.5cm,clip,width=\textwidth]{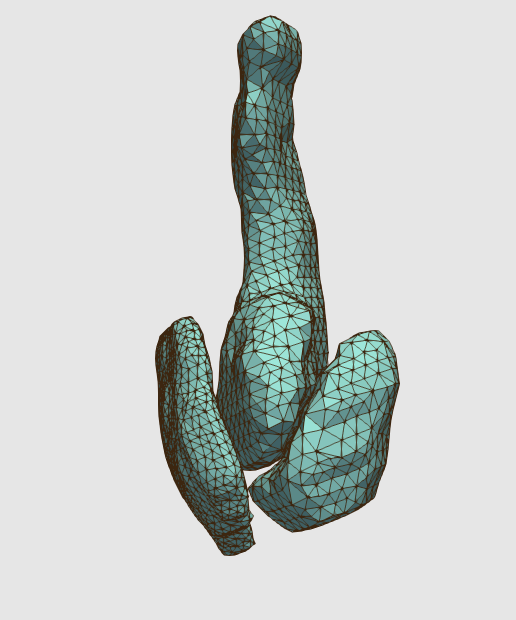}
\caption*{(a)}
\end{subfigure}
\begin{subfigure}[t]{0.4\textwidth}
\centering
\includegraphics[trim=4cm 2cm 3cm 0.5cm,clip,width=\textwidth]{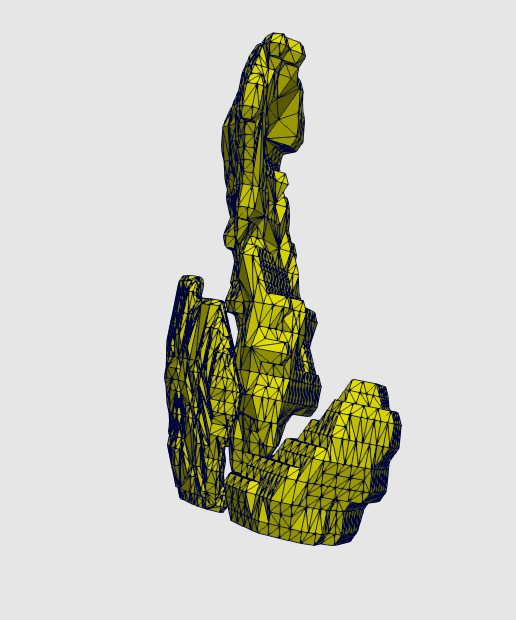}
\caption*{(b)}
\end{subfigure}
\caption{\label{fig:biocardMulti} (a) Template and (b) target, each of them the union of three surfaces (left: entorhinal cortex, center: hippocampus, right: amygdala)} 
\end{figure}

\begin{figure}
\centering
\begin{subfigure}[t]{0.22\textwidth}
\centering
\includegraphics[trim=4cm 2cm 3cm 0.5cm,clip,width=\textwidth]{biocardMultiTemplate.png}\
\caption*{$t=0$}
\end{subfigure}
\begin{subfigure}[t]{0.22\textwidth}
\centering
\includegraphics[trim=4cm 2cm 3cm 0.5cm,clip,width=\textwidth]{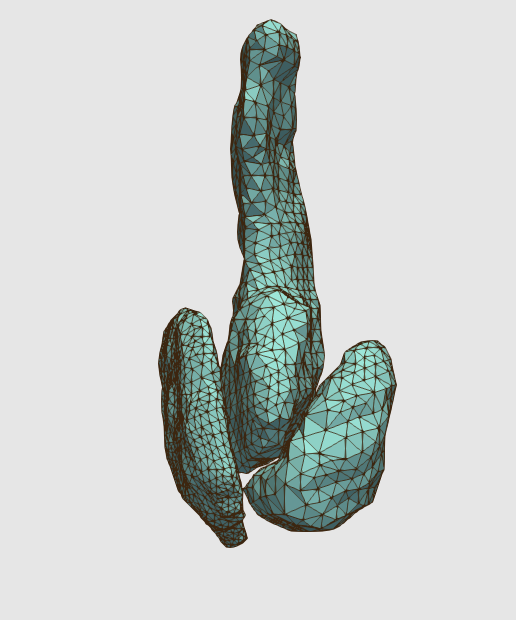}\
\caption*{$t=0.3$}
\end{subfigure}
\begin{subfigure}[t]{0.22\textwidth}
\centering
\includegraphics[trim=4cm 2cm 3cm 0.5cm,clip,width=\textwidth]{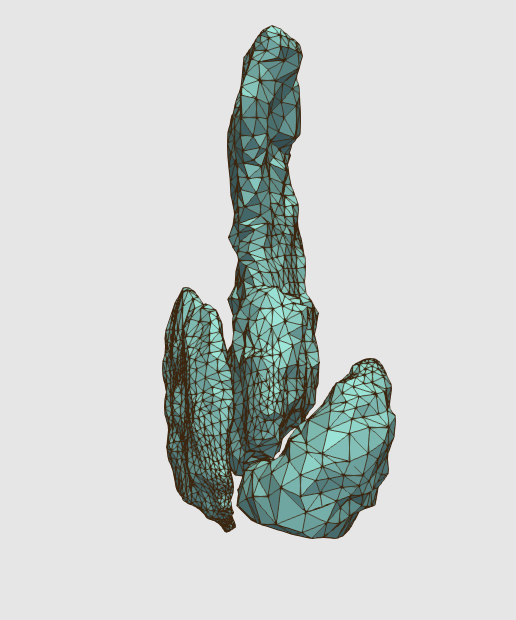}\
\caption*{$t=0.7$}
\end{subfigure}
\begin{subfigure}[t]{0.22\textwidth}
\centering
\includegraphics[trim=4cm 2cm 3cm 0.5cm,clip,width=\textwidth]{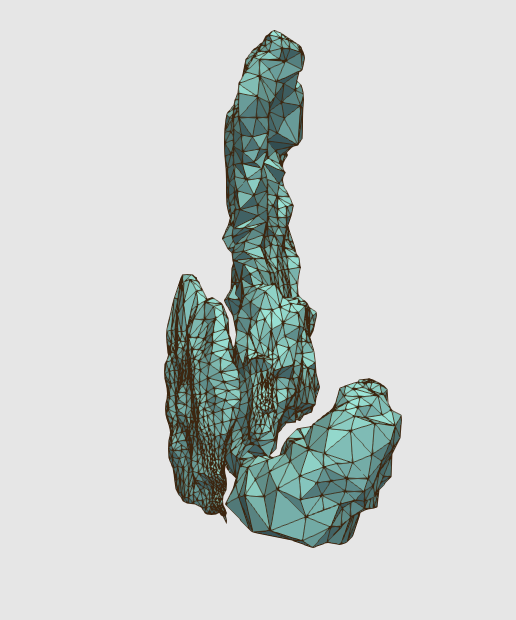}\\
\caption*{$t=1.0$}
\end{subfigure}
\begin{subfigure}[t]{0.22\textwidth}
\centering
\includegraphics[trim=4cm 2cm 3cm 0.5cm,clip,width=\textwidth]{biocardMultiTemplate.png}\
\caption*{$t=0$}
\end{subfigure}
\begin{subfigure}[t]{0.22\textwidth}
\centering
\includegraphics[trim=4cm 2cm 3cm 0.5cm,clip,width=\textwidth]{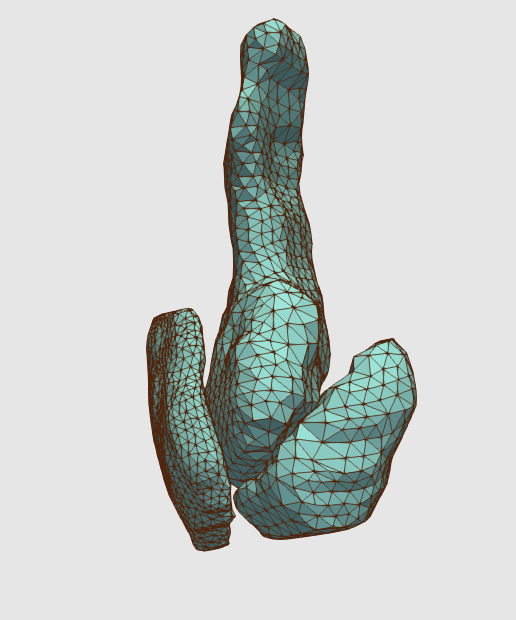}\
\caption*{$t=0.3$}
\end{subfigure}
\begin{subfigure}[t]{0.22\textwidth}
\centering
\includegraphics[trim=4cm 2cm 3cm 0.5cm,clip,width=\textwidth]{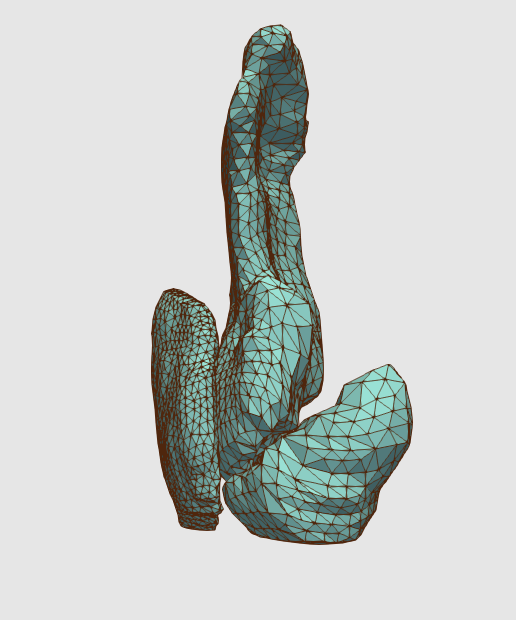}\
\caption*{$t=0.7$}
\end{subfigure}
\begin{subfigure}[t]{0.22\textwidth}
\centering
\includegraphics[trim=4cm 2cm 3cm 0.5cm,clip,width=\textwidth]{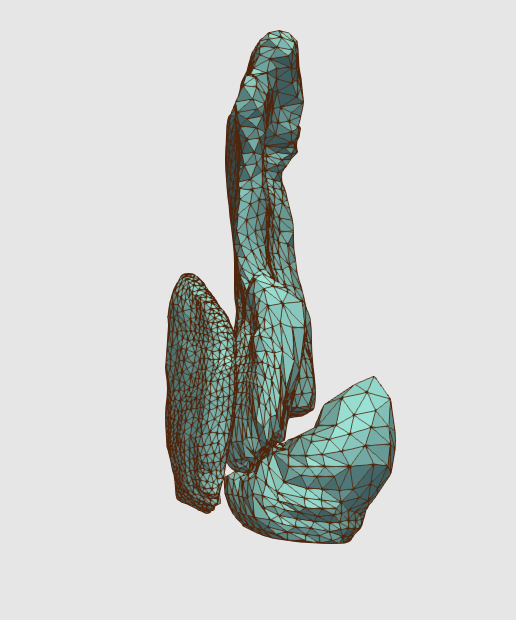}\\
\caption*{$t=1.0$}
\end{subfigure}
\begin{subfigure}[t]{0.22\textwidth}
\centering
\includegraphics[trim=4cm 2cm 3cm 0.5cm,clip,width=\textwidth]{biocardMultiTemplate.png}\
\caption*{$t=0$}
\end{subfigure}
\begin{subfigure}[t]{0.22\textwidth}
\centering
\includegraphics[trim=4cm 2cm 3cm 0.5cm,clip,width=\textwidth]{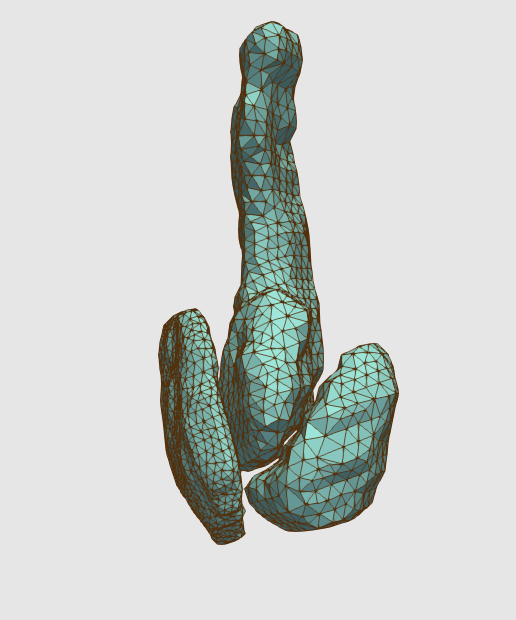}\
\caption*{$t=0.3$}
\end{subfigure}
\begin{subfigure}[t]{0.22\textwidth}
\centering
\includegraphics[trim=4cm 2cm 3cm 0.5cm,clip,width=\textwidth]{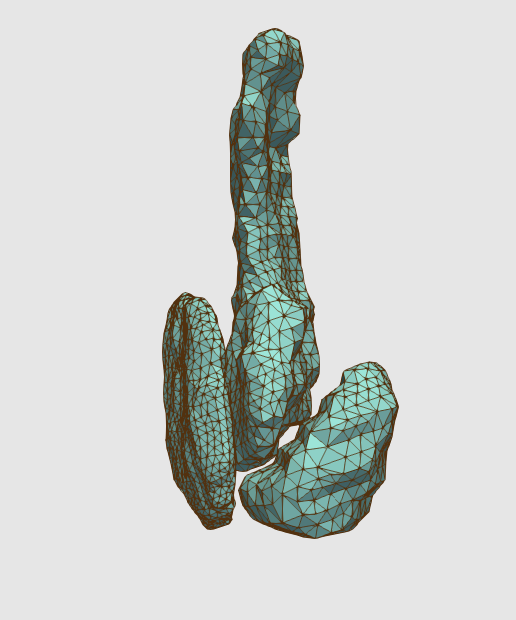}\
\caption*{$t=0.7$}
\end{subfigure}
\begin{subfigure}[t]{0.22\textwidth}
\centering
\includegraphics[trim=4cm 2cm 3cm 0.5cm,clip,width=\textwidth]{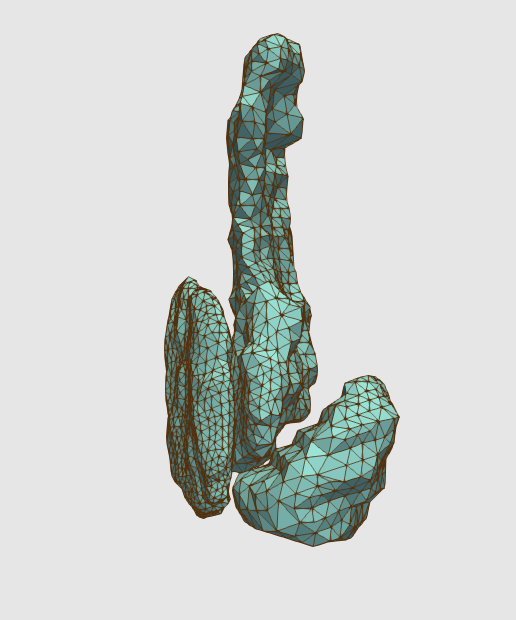}\\
\caption*{$t=1.0$}
\end{subfigure}
\caption{\label{fig:biocardMulti.geod} Line segments: geodesics. First row: LDDMM. Second row: Hybrid LDDMM.}
\end{figure}

Our example uses the same data as the one presented in Figures 8 and 9 of \cite{arguillere2016registration}. It includes three shapes (see Figure \ref{fig:biocardMulti}) who are relatively close to each other (the hippocampus and the amygdala are actually slightly overlapping in the target). If one uses standard LDDMM with a small kernel ($a \simeq d/45$, where $d$ is the height of the hippocampus), as illustrated in the first row of Figure \ref{fig:biocardMulti.geod}, the diffeomorphism has undesirable properties, crunching parts of the surfaces (such as the front of the hippocampus, the bottom of the entorhinal cortex ---which even has a residual spike--- and the top of the amygdala) to match the target. With a larger kernel width ($a\simeq d/6$) the three shapes are transformed as if they formed one single object, resulting in large reparametrization of the surfaces when they move along each other, because points that were nearby initially tend to have similar motions. This is illustrated in the second row of Figure \ref{fig:biocardMulti.geod}. The third row provides the geodesic obtained with the hybrid norm, with the same small kernel width as in the first row, but with \eqref{eq:h1.surf} penalizing large deformations on the surfaces. In this case, the surfaces move nicely along each other, without requiring large reparametrization, except those required by the change in their respective shapes. 

\section{Conclusion}

Our results illustrate several advantages of combining the standard LDDMM approach with geometrically inspired norms on spaces of curves and surfaces. This very simple concept allows for much more modeling accuracy and flexibility, with a moderate computational impact. This is especially useful when dealing with complex configurations of shapes, as we saw in our examples. 

There is clearly still room for future work and development, including the use of higher-order norms for curves and surfaces,  and guidelines on which norm to use in specific applications. A version of the method for image matching is another important direction to be explored in the future. Formally, this requires defining $\normi{h}_q$, where both $q$ and $h$ are scalar functions in $\mR^d$ and using the extra term $\normi{v^T\nabla q}_q$ in the Riemannian norm. One option worth exploring is
\[
\normi{h}_q^2 = \int_{\mR^d} |\nabla h(x)|^2 w_q(x)\,dx
\]
where $w_q$ is a ``weight function'' that depends on $q$, making, for example, deformations more costly in gray/white matter regions than within cerebro-spinal fluid in brain mapping. This will be addressed in future work.

\bibliography{refShape}

\begin{thebibliography}{10}

\bibitem{arguillere2015abstract}
Sylvain Arguill{\`e}re.
\newblock The abstract setting for shape deformation analysis and lddmm
  methods.
\newblock In Frank Nielsen and Fr{\'e}d{\'e}ric Barbaresco, editors, {\em
  Geometric Science of Information: Second International Conference, GSI 2015,
  Palaiseau, France, October 28-30, 2015, Proceedings}, pages 159--167, Cham,
  2015. Springer International Publishing.

\bibitem{arguillere2014shape}
Sylvain Arguill{\`e}re, Emmanuel Tr{\'e}lat, Alain Trouv{\'e}, and Laurent
  Younes.
\newblock Shape deformation analysis from the optimal control viewpoint.
\newblock {\em Journal de Math{\'e}matiques Pures et Appliqu{\'e}es},
  104(1):139--178, 2015.

\bibitem{arguillere2016registration}
Sylvain Arguill{\`e}re, Emmanuel Tr{\'e}lat, Alain Trouv{\'e}, and Laurent
  Younes.
\newblock Registration of multiple shapes using constrained optimal control.
\newblock {\em SIAM Journal on Imaging Sciences}, 9(1):344--385, 2016.

\bibitem{af00}
J.~Ashburner and K.~J. Friston.
\newblock Voxel based morphometry -- the methods.
\newblock {\em Neuroimage}, 11(6):805--821, 2000.

\bibitem{ashburner2007fast}
John Ashburner.
\newblock A fast diffeomorphic image registration algorithm.
\newblock {\em Neuroimage}, 38(1):95--113, 2007.

\bibitem{ashburner2011diffeomorphic}
John Ashburner and Karl~J Friston.
\newblock Diffeomorphic registration using geodesic shooting and gauss--newton
  optimisation.
\newblock {\em NeuroImage}, 55(3):954--967, 2011.

\bibitem{avants2004geodesic}
Brian Avants and James~C Gee.
\newblock Geodesic estimation for large deformation anatomical shape averaging
  and interpolation.
\newblock {\em Neuroimage}, 23:S139--S150, 2004.

\bibitem{bauer2011new}
Martin Bauer and Martins Bruveris.
\newblock A new {Riemannian} setting for surface registration setting for
  surface registration.
\newblock In {\em Proceedings of the Third International Workshop on
  Mathematical Foundations of Computational Anatomy-Geometrical and Statistical
  Methods for Modelling Biological Shape Variability}, pages 182--193, 2011.

\bibitem{bauer2014constructing}
Martin Bauer, Martins Bruveris, Stephen Marsland, and Peter~W Michor.
\newblock Constructing reparameterization invariant metrics on spaces of plane
  curves.
\newblock {\em Differential Geometry and its Applications}, 34:139--165, 2014.

\bibitem{bauer2014overview}
Martin Bauer, Martins Bruveris, and Peter~W Michor.
\newblock Overview of the geometries of shape spaces and diffeomorphism groups.
\newblock {\em Journal of Mathematical Imaging and Vision}, 50(1-2):60--97,
  2014.

\bibitem{bauer2011sobolev}
Martin Bauer, Philipp Harms, and Peter~W Michor.
\newblock Sobolev metrics on shape space of surfaces.
\newblock {\em Journal of Geometric Mechanics (2011), 389-438}, pages 389--438,
  2011.

\bibitem{beg2005computing}
M~Faisal Beg, Michael~I Miller, Alain Trouv{\'e}, and Laurent Younes.
\newblock Computing large deformation metric mappings via geodesic flows of
  diffeomorphisms.
\newblock {\em International journal of computer vision}, 61(2):139--157, 2005.

\bibitem{cao2005large}
Yan Cao, Michael~I Miller, Raimond~L Winslow, and Laurent Younes.
\newblock Large deformation diffeomorphic metric mapping of vector fields.
\newblock {\em Medical Imaging, IEEE Transactions on}, 24(9):1216--1230, 2005.

\bibitem{ceritoglu2009multi}
Can Ceritoglu, Kenichi Oishi, Xin Li, Ming-Chung Chou, Laurent Younes, Marilyn
  Albert, Constantine Lyketsos, Peter van Zijl, Michael~I Miller, and Susumu
  Mori.
\newblock Multi-contrast large deformation diffeomorphic metric mapping for
  diffusion tensor imaging.
\newblock {\em Neuroimage}, 47(2):618--627, 2009.

\bibitem{charon2013varifold}
N.~Charon and A.~Trouv{\'e}.
\newblock The varifold representation of nonoriented shapes for diffeomorphic
  registration.
\newblock {\em SIAM Journal on Imaging Sciences}, 6(4):2547--2580, 2013.

\bibitem{christensen1996deformable}
Gary~E. Christensen, Richard~D. Rabbitt, and Michael~I. Miller.
\newblock Deformable templates using large deformation kinematics.
\newblock {\em Image Processing, IEEE Transactions on}, 5(10):1435--1447, 1996.

\bibitem{droske2004variational}
Marc Droske and Martin Rumpf.
\newblock A variational approach to nonrigid morphological image registration.
\newblock {\em SIAM Journal on Applied Mathematics}, 64(2):668--687, 2004.

\bibitem{glaunes2008large}
Joan Glaun{\`e}s, Anqi Qiu, Michael~I Miller, and Laurent Younes.
\newblock Large deformation diffeomorphic metric curve mapping.
\newblock {\em International journal of computer vision}, 80(3):317--336, 2008.

\bibitem{hernandez2008comparing}
Monica Hernandez, Salvador Olmos, and Xavier Pennec.
\newblock Comparing algorithms for diffeomorphic registration: Stationary lddmm
  and diffeomorphic demons.
\newblock In {\em 2nd MICCAI Workshop on Mathematical Foundations of
  Computational Anatomy}, pages 24--35, 2008.

\bibitem{holm1998euler}
Darryl~D Holm, Jerrold~E Marsden, and Tudor~S Ratiu.
\newblock The {Euler--Poincar{\'e}} equations and semidirect products with
  applications to continuum theories.
\newblock {\em Advances in Mathematics}, 137(1):1--81, 1998.

\bibitem{holm2009euler}
Darryl~D Holm, Alain Trouv{\'e}, and Laurent Younes.
\newblock The {Euler-Poincar{\'e}} theory of metamorphosis.
\newblock {\em Quarterly of Applied Mathematics}, 97:661--685, 2009.

\bibitem{jos97}
S.~Joshi.
\newblock {\em Large Deformation Diffeomorphisms and {Gaussian} Random Fields
  for Statistical Characterization of Brain Sub-manifolds}.
\newblock PhD thesis, Sever institute of technology, Washington University,
  1997.

\bibitem{joshi2000Landmark}
Sarang~C. Joshi and Michael~I. Miller.
\newblock {Landmark matching via large deformation diffeomorphisms}.
\newblock {\em IEEE Transactions on Image Processing}, 9:1357--1370, 2000.

\bibitem{ksmj04}
E.~Klassen, A.~Srivastava, W.~Mio, and S.~H. Joshi.
\newblock Analysis of planar shapes using geodesic paths on shape spaces.
\newblock {\em IEEE Trans. Pattern Anal. Mach. Intell.}, 26(3):372--383, 2004.

\bibitem{klein2009evaluation}
Arno Klein, Jesper Andersson, Babak~A. Ardekani, John Ashburner, Brian Avants,
  Ming-Chang Chiang, Gary~E. Christensen, D.~Louis Collins, James Gee, Pierre
  Hellier, Joo~Hyun Song, Mark Jenkinson, Claude Lepage, Daniel Rueckert, Paul
  Thompson, Tom Vercauteren, Roger~P. Woods, J.~John Mann, and Ramin~V. Parsey.
\newblock Evaluation of 14 nonlinear deformation algorithms applied to human
  brain \{MRI\} registration.
\newblock {\em NeuroImage}, 46(3):786 -- 802, 2009.

\bibitem{krall1949new}
Hi~L Krall and Orrin Frink.
\newblock A new class of orthogonal polynomials: The bessel polynomials.
\newblock {\em Transactions of the American Mathematical Society},
  65(1):100--115, 1949.

\bibitem{kurtek2011parameterization}
Sebastian Kurtek, Eric Klassen, Zhaohua Ding, Sandra~W Jacobson, Joseph~L
  Jacobson, Malcolm~J Avison, and Anuj Srivastava.
\newblock Parameterization-invariant shape comparisons of anatomical surfaces.
\newblock {\em IEEE Transactions on Medical Imaging}, 30(3):849--858, 2011.

\bibitem{mm06a}
P.~W. Michor and D.~Mumford.
\newblock Riemannian geometries on spaces of plane curves.
\newblock {\em J. Eur. Math. Soc.}, 8:1--48, 2006.

\bibitem{mm07}
P.~W. Michor and D.~Mumford.
\newblock An overview of the {Riemannian} metrics on spaces of curves using the
  {Hamiltonian} approach.
\newblock {\em Applied and Computational Harmonic Analysis}, 23(1):74--113,
  2007.

\bibitem{miller2006geodesic}
Michael~I Miller, Alain Trouv{\'e}, and Laurent Younes.
\newblock Geodesic shooting for computational anatomy.
\newblock {\em Journal of mathematical imaging and vision}, 24(2):209--228,
  2006.

\bibitem{miller2001group}
Michael~I. Miller and Laurent Younes.
\newblock Group actions, homeomorphisms, and matching: A general framework.
\newblock {\em International Journal of Computer Vision}, 41(1-2):61--84, 2001.

\bibitem{pantazis2004statistical}
Dimitrios Pantazis, Richard~M Leahy, Thomas~E Nichols, and Martin Styner.
\newblock Statistical surface-based morphometry using a nonparametric approach.
\newblock In {\em Biomedical Imaging: Nano to Macro, 2004. IEEE International
  Symposium on}, pages 1283--1286. IEEE, 2004.

\bibitem{risser2011simultaneous}
Laurent Risser, F~Vialard, Robin Wolz, Maria Murgasova, Darryl~D Holm, and
  Daniel Rueckert.
\newblock Simultaneous multi-scale registration using large deformation
  diffeomorphic metric mapping.
\newblock {\em Medical Imaging, IEEE Transactions on}, 30(10):1746--1759, 2011.

\bibitem{styner2006framework}
Martin Styner, Ipek Oguz, Shun Xu, Christian Brechb{\"u}hler, Dimitrios
  Pantazis, James~J Levitt, Martha~E Shenton, and Guido Gerig.
\newblock Framework for the statistical shape analysis of brain structures
  using spharm-pdm.
\newblock {\em The insight journal}, 1071:242, 2006.

\bibitem{trouve2005metamorphoses}
Alain Trouv{\'e} and Laurent Younes.
\newblock Metamorphoses through {Lie} group action.
\newblock {\em Foundations of Computational Mathematics}, 5(2):173--198, 2005.

\bibitem{vercauteren2009diffeomorphic}
Tom Vercauteren, Xavier Pennec, Aymeric Perchant, and Nicholas Ayache.
\newblock Diffeomorphic demons: Efficient non-parametric image registration.
\newblock {\em NeuroImage}, 45(1):S61--S72, 2009.

\bibitem{yezzi2004metrics}
Anthony Yezzi and Andrea Mennucci.
\newblock Metrics in the space of curves.
\newblock {\em arXiv preprint math/0412454}, 2004.

\bibitem{ymsm07}
L.~Younes, P.~Michor, J.~Shah, and D.~Mumford.
\newblock A metric on shape spaces with explicit geodesics.
\newblock {\em Rend. Lincei Mat. Appl.}, 9:25--57, 2008.

\bibitem{younes1998computable}
Laurent Younes.
\newblock Computable elastic distances between shapes.
\newblock {\em SIAM Journal on Applied Mathematics}, 58(2):565--586, 1998.

\bibitem{younes2007jacobi}
Laurent Younes.
\newblock Jacobi fields in groups of diffeomorphisms and applications.
\newblock {\em Quarterly of applied mathematics}, 65(1):113--134, 2007.

\bibitem{younes2009evolutions}
Laurent Younes, Felipe Arrate, and Michael~I Miller.
\newblock Evolutions equations in computational anatomy.
\newblock {\em NeuroImage}, 45(1):S40--S50, 2009.

\end{thebibliography}
\bibliographystyle{plain}

\end{document}